\documentclass[11pt]{article}
\usepackage[T1]{fontenc}
\usepackage{lmodern,amsmath,amssymb,amsthm,mathtools}
\usepackage[margin=1in]{geometry}
\usepackage{booktabs,array,longtable,enumitem,microtype}
\usepackage[hidelinks]{hyperref}
\usepackage{xurl}
\newtheorem{theorem}{Theorem}[section]
\newtheorem{proposition}[theorem]{Proposition}
\newtheorem{lemma}[theorem]{Lemma}

\theoremstyle{definition}
\newtheorem{definition}[theorem]{Definition}
\newtheorem{assumption}[theorem]{Input}
\newtheorem*{example*}{Example}
\theoremstyle{remark}
\newtheorem{remark}[theorem]{Remark}
\numberwithin{equation}{section}
\newcommand{\HH}{\mathbb H}
\newcommand{\R}{\mathbb R}
\newcommand{\Q}{\mathbb Q}

\newcommand{\ip}[2]{\langle #1,#2\rangle}
\newcommand{\norm}[1]{\lVert #1\rVert}

\newcommand{\one}{\mathbf 1}
\title{Linear programming certificates for dimension bounds\\
of hyperbolic Coxeter polytopes}
\author{Alexander Kolpakov}
\date{19 November 2026}
\begin{document}
\maketitle

\begin{abstract}
We use linear programming to derive upper bounds for the dimension of
hyperbolic Coxeter polytopes from local geometric and arithmetic
estimates. Rational certificates give the bounds $57$ for finite-volume
polytopes, $17$ for compact polytopes, and $15$ for compact arithmetic
polytopes. In the finite-volume case, the certificate is a weight
function on incidences between vertices and two-dimensional faces.
In the compact cases, the inequalities concern squared norms of
orthogonal projections of a fixed tangent vector at a vertex
farthest from an interior point.
We give a common formulation of these linear programs, prove
monotonicity in the dimension threshold, and describe the exact
verification of the certificates. The local estimates are taken from
the companion geometric proofs.
\end{abstract}

\tableofcontents

\section{Introduction}
\label{sec:introduction}

Upper bounds for the dimension of a hyperbolic Coxeter polytope can be
proved by comparing local inequalities with a global incidence count
or a bound on a vector norm. We formulate these comparisons as linear
programs. Their constraints are necessary conditions for a polytope
to exist, and their infeasibility proves a dimension bound. The dual
linear programs describe the combinations of local inequalities used
in the proof.

We recover the following bounds from the local results of the
companion paper and supplement~\cite{BoundsArticle,BoundsSupplement}:
\begin{equation}\label{eq:three-bounds}
 \begin{aligned}
 n&\le57&&\text{in finite volume},\\
 n&\le17&&\text{in the compact case},\\
 n&\le15&&\text{in the compact arithmetic case}.
 \end{aligned}
\end{equation}
Guo and Zheng~\cite{GuoZheng2026} prove the general finite-volume
bound $n\le81$. Our purpose here is to give a uniform linear
programming treatment of~\eqref{eq:three-bounds}.

Throughout the paper, a polytope is convex, finite-sided and
full-dimensional. Hyperbolic space has curvature $-1$. A Coxeter
polytope has all dihedral angles of the form $\pi/m$, where $m\ge2$
is an integer. Its reflection group is generated by reflections in
the supporting hyperplanes of its facets. Arithmeticity always refers
to this group. We write Supplements FV, C and A for Parts I/IV, II/V
and III/VI, respectively, of~\cite{BoundsSupplement}. The local
results used below are stated with their hypotheses; their geometric
proofs are given in those supplements.

\begin{theorem}[Linear programming proofs of the dimension bounds]
\label{thm:recovery}
The local estimates and classifications stated below give the following
rational certificates.
\begin{enumerate}[label=\textup{(\roman*)}]
\item For every integer $n\ge58$, there is a nonnegative weight
function on vertex--face incidences satisfying the face and vertex
inequalities of Proposition~\ref{fv:charge-certificate}.
Consequently the finite-volume incidence problem is infeasible.
\item The normalized linear program for the squared projection norms
of compact polytopes, defined in~\eqref{eq:budget-deficit}, is
infeasible on the exhaustive family of configurations in dimensions
$18,\ldots,29$. Its dual has a rational solution proving this
infeasibility.
\item The same conclusion holds for compact arithmetic polytopes in
dimensions $16,\ldots,29$. At dimensions $16$ and $17$ the program
uses the arithmetic classifications and inequalities stated below;
at dimensions $18,\ldots,29$ it uses the compact configurations
and inequalities from part \textup{(ii)}.
\end{enumerate}
Together with Vinberg's compact bound $n\le29$, these certificates
imply~\eqref{eq:three-bounds}.
\end{theorem}

Theorem~\ref{thm:homogeneous-certificate} gives the common linear
argument. A hypothetical polytope determines nonnegative variables
satisfying homogeneous local constraints and a strict global
inequality. Dividing by its positive excess gives a solution of a
normalized feasibility problem. A dual inequality contradicts that
normalization. Proposition~\ref{prop:threshold} proves monotonicity:
raising the lower bound on the dimension removes configurations
from one fixed family.

\begin{center}\small
\begin{tabular}{>{\raggedright\arraybackslash}p{.15\linewidth}
 >{\raggedright\arraybackslash}p{.41\linewidth}
 >{\raggedright\arraybackslash}p{.34\linewidth}}
\toprule
Case&Quantities in the linear program&Dual inequalities\\\midrule
Finite volume&The same corner functions summed over faces and over
 vertices&Lower bounds for face sums and upper bounds for vertex sums\\
Compact&Squared norms of orthogonal projections of a fixed vector,
 whose total is less than one&Nonnegative combinations of local
 inequalities with coefficient at most one on each variable\\
Compact arithmetic&The same projection norms, subject to additional
 arithmetic restrictions&The same dual program, with the arithmetic
 inequalities and excluded configurations\\\bottomrule
\end{tabular}
\end{center}

The local estimates and the exhaustive lists of configurations are
hypotheses of the linear argument. This includes the reductions
leading to the two arithmetic lists. Conditions on the searched dual
coefficients have a different role: they ensure that a linear
combination proves the desired contradiction. We distinguish these
conditions from geometric restrictions on a polytope wherever both
occur in the same formulation.

Each local estimate is used over its stated range of ranks, labels
and other components. Ordinary dihedral components retain their full
integer parameter range. Configurations satisfying these necessary
conditions need not be geometrically realizable; every polytope must
nevertheless give one of them. All terminating decimals in this paper
denote exact rational numbers.

\section{Geometric notation and linear programming}
\label{sec:framework}

\subsection{Facet normals and Coxeter diagrams}

We use a Lorentz space of signature $(-,+,\ldots,+)$, with hyperbolic
space the future unit hyperboloid. The supporting hyperplane of a
facet will be called a wall. It has an outward spacelike unit
normal $u$ with $\ip{u}{u}=1$, and the chamber inequality is
$\ip{y}{u}\le0$. For two facets meeting at angle $\pi/m$, the corresponding Gram
entry is $-\cos(\pi/m)$. At a finite vertex the Gram matrix of the
incident facets is positive definite. The spans of the connected
components of its Coxeter diagram are mutually orthogonal and together
form the Euclidean tangent space.

A spherical Coxeter diagram has positive definite Gram matrix. A
connected affine diagram has positive semidefinite Gram matrix of
corank one. We use the standard list of irreducible spherical diagrams:
\[
 A_r\ (r\ge1),\quad B_r\ (r\ge3),\quad D_r\ (r\ge4),\quad
 E_6,E_7,E_8,F_4,H_3,H_4,\quad I_2(m)\ (m\ge4).
\]
Here $A_2=I_2(3)$; the isomorphic diagrams $B_2$ and $D_3$ are not
listed separately. The \emph{order}
of a diagram is its number of nodes, and its \emph{rank} is the rank
of its Gram matrix. They agree for a spherical diagram. A connected
affine diagram with $r+1$ nodes has Gram rank $r$. The simplex factor
of an ideal vertex link that it determines has dimension $r$.
We use this distinction between order and rank in the vertex-figure
formulas below.

All facets and facet normals appearing in a local geometric
statement belong to the polytope $P$, unless explicitly introduced
as auxiliary objects. A facet opposite a vertex along a specified
edge means a facet containing the other endpoint but not that
vertex. Replacing a normal by a projected vector does not produce
another facet of $P$. An auxiliary point asserted to lie in $P$
must satisfy every defining facet inequality. The Coxeter and arithmetic hypotheses
are imposed on $P$; they are not assumed for an arbitrary face.

\subsection{A vector at a farthest vertex}

Let $P$ be compact. Fix a point $z\in\operatorname{int}P$ and choose
a vertex $p$ at maximal distance from $z$. Put
$H=-\ip{z}{p}=\cosh d(z,p)$.
The decomposition
\begin{equation}\label{eq:original-vector}
                   z=H(p-x),\qquad x\in p^\perp
\end{equation}
defines a vector $x$ in the tangent space at $p$. Indeed,
$x=p-z/H$ is orthogonal to $p$, and $z^2=p^2=-1$ gives
\begin{equation}\label{eq:original-budget}
                         \norm{x}^2=1-H^{-2}<1.
\end{equation}
The Lorentz form restricts to a Euclidean inner product on $p^\perp$.
Let $W_1,\ldots,W_k$ be the spans of the connected components of the
Coxeter diagram at $p$, and define
\begin{equation}\label{eq:component-energies}
 e_j=\norm{\operatorname{proj}_{W_j}x}^2,\qquad
 e_j\ge0,\qquad \sum_{j=1}^ke_j=\norm{x}^2<1.
\end{equation}
Thus $e_j$ is a squared norm, rather than a projection vector. A sum
of several $e_j$ is the squared norm of the projection onto the
corresponding orthogonal sum of subspaces. Throughout both compact
arguments, all these quantities refer to the fixed vector $x$.
Auxiliary points used in a local proof do not supply further terms
in~\eqref{eq:component-energies}.

\subsection{Classes of configurations and local inequalities}

\begin{definition}\label{def:configuration-family}
An \emph{exhaustive family of configuration classes} in dimension $n$
is a finite collection of sets $\mathcal D_b$, with $n_b=n$, whose
union contains every geometric configuration under consideration.
Each class is specified by all the hypotheses of the local results
used on it. A class may include a whole parameter range when its
inequalities hold uniformly throughout that range.
\end{definition}

The classes may overlap. A geometric configuration can be assigned
to any class containing it. Some formal configurations in the classes
may have no geometric realization; their inclusion preserves the
necessary direction of the argument.

A local norm inequality specifies a set $S$ of the orthogonal
components and a rational number $c\ge0$ such that
\begin{equation}\label{eq:supported-row}
                            \sum_{j\in S}e_j\ge c
\end{equation}
on a stated configuration class. We call $S$ the \emph{support} of
the inequality. A lower bound for the squared norm of a projection onto a subspace of one
component also bounds the squared norm on that component, by
nonnegativity of the omitted terms. Bounds on overlapping subspaces
must still be combined with their repeated coefficients accounted for.

A local theorem may instead show that a specified configuration
cannot occur. We call this a \emph{local exclusion}. Its geometric
and arithmetic hypotheses form part of the definition of the class.
In the homogeneous formulation its scaling variable is set to zero.
Such an exclusion is used only on the stated class, whereas exclusion
of every class in a given dimension is the conclusion of the LP.

\subsection{Homogeneous linear programs}

We first state the finite-dimensional alternative underlying the
linear certificates.

\begin{lemma}\label{lem:cone-alternative}
Let $M$ be a rational matrix and $d$ a rational column vector of
compatible size. Exactly one of the following
systems has a solution:
\begin{align}
 Mu=d,\quad u\ge0; &\label{eq:cone-primal}\\
 M^Ty\ge0,\quad d^Ty=-1. &\label{eq:cone-dual}
\end{align}
When either system is feasible, it has a rational solution.
\end{lemma}
\begin{proof}
Let $K$ be the cone generated by the columns of $M$. It is a closed
polyhedral cone. If $d\notin K$, separation gives a vector $y$ with
$y^Tk\ge0$ on $K$ and $y^Td<0$; rescaling gives
\eqref{eq:cone-dual}. If $d\in K$, its nonnegative column
coefficients give~\eqref{eq:cone-primal}. Both systems cannot be
feasible, since they would give
$-1=d^Ty=u^TM^Ty\ge0$.

For rationality, a nonempty polyhedron defined over $\Q$ has a
rational point: restrict to the affine hull of a minimal face,
solve its rational equality system, and approximate a relative
interior point by rational points in that affine space. All remaining
strict inequalities persist under a sufficiently small approximation.
Apply this to the two displayed rational systems.
\end{proof}

Replacing a free variable by a difference of nonnegative variables
and introducing slack variables puts a finite linear system in this
form. Below we also write its dual inequalities in terms of the
geometric quantities they estimate.

Let $W$ be a real vector space and $I$ an index set. To each $i\in I$
assign an integer $n_i$, a polyhedral cone $C_i\subset\R^{k_i}$ of local variables,
a linear balance map $T_i:\R^{k_i}\longrightarrow W$, and a linear
functional $d_i:\R^{k_i}\longrightarrow\R$. The label $n_i$ denotes
the dimension of the polytope, even when $i$ describes a local
configuration. The cone $C_i$ is defined by homogeneous necessary
conditions. The
map $T_i$ records the terms in the linear identities to be imposed,
and $d_i$ their contribution to a strict inequality. Fix these data
independently of the lower dimension threshold.

For an integer $N$ consider
\begin{equation}\label{eq:common-feasibility}
 \begin{gathered}
 u_i\in C_i\quad(n_i\ge N),\qquad
 \sum_{n_i\ge N}T_i u_i=0,\qquad
 \sum_{n_i\ge N}d_i(u_i)=1.
 \end{gathered}
\end{equation}
Only finitely many of the vectors $u_i$ may be nonzero. For infinitely many
dimensions, take $W$ to be a direct sum with separate coordinates for
each dimension. Incidence identities from different dimensions then
remain separate. A finite-sided polytope contributes only finitely
many terms.

\begin{theorem}[Normalization and dual certificate]
\label{thm:homogeneous-certificate}
Suppose every polytope under consideration of dimension $n\ge N$
determines a finitely supported family of vectors $v_i\in C_i$,
zero whenever $n_i\ne n$, with
\[
             \sum_iT_iv_i=0,\qquad D:=\sum_id_i(v_i)>0.
\]
Then~\eqref{eq:common-feasibility} has a solution.
A linear functional $\lambda$ on $W$ satisfying
\begin{equation}\label{eq:common-certificate}
       \lambda(T_iu)-d_i(u)\ge0
       \qquad(u\in C_i,\ n_i\ge N)
\end{equation}
proves infeasibility and hence excludes all these polytopes.
If $I$ is finite, $W$ is finite-dimensional and all the cones, maps
and functionals are rational, infeasibility is equivalent to the
existence of such a rational certificate.
\end{theorem}
\begin{proof}
Set $u_i=v_i/D$ to obtain the normalized solution. If a feasible
point and a dual certificate both existed, summing
\eqref{eq:common-certificate} at that point would give
$0-1\ge0$, a contradiction.
For the finite converse, express each rational polyhedral cone by
finitely many rational generators and apply
Lemma~\ref{lem:cone-alternative} to their images
$(T_iu,d_i(u))$, with target $(0,1)$. Its separating functional has
the form $(\lambda,-1)$, exactly as in
\eqref{eq:common-certificate}.
\end{proof}

\begin{proposition}[Monotonicity in the dimension threshold]\label{prop:threshold}
In the fixed formulation~\eqref{eq:common-feasibility}, feasibility
at $N'\ge N$ implies feasibility at $N$. Thus infeasibility at $N$
implies infeasibility at every larger threshold. A certificate at
$N$ remains a certificate after restricting to $n_i\ge N'$.
\end{proposition}
\begin{proof}
Extend a feasible point by zeros on the removed indices. Every cone
contains zero, and both sums are unchanged. The same dual inequalities hold after restriction to the remaining
indices.
\end{proof}

For an infinite index set only the inequality on finite sums is
needed; the converse duality statement is restricted to finite
systems. If the configuration classes cover a bounded range of
dimensions, a separate theorem must exclude the higher dimensions.
Section~\ref{sec:thresholds} states this additional information for
each application. No passage to a lower-dimensional Coxeter face
is involved.

\section{Finite volume: weights on vertex--face incidences}
\label{fv:section}

We assign a weight to each incidence of a vertex and a two-dimensional
face. Summing these weights first over faces and then over vertices
gives the comparison used to bound the dimension. We formulate the
required lower bounds at faces and upper bounds at vertices as a
linear program. A weight function satisfying these inequalities
contradicts the incidence inequality stated below. We also minimize the required vertex bound
within a specified family of weight functions.

\subsection{Vertex figures and the incidence inequality}
\label{fv:incidence}

Fix an integer $n\ge58$. A \emph{planar angle}, or \emph{corner},
is an incidence $\kappa=(v,F)$ of a vertex $v$ and a two-dimensional
face $F$.
The number of sides of $F$ is denoted by $s(F)$, and
\begin{equation}\label{fv:defect}
 D_5(P)=\sum_{\dim F=2}(5-s(F)).
\end{equation}
A triangle contributes two, a quadrilateral contributes one, and a
face with at least five sides contributes a nonpositive number.

We use the standard vertex-figure description of a finite-volume
Coxeter polytope. At a finite vertex the figure is a simplex. At an
ideal vertex a horospherical figure is a product
\begin{equation}\label{fv:link}
 L_v=\prod_{i=1}^k\Delta^{r_i},\qquad
 r_i\ge1,\qquad \sum_i r_i=n-1.
\end{equation}
Here $r_i$ is the dimension of a simplex factor, not the number of
nodes of its affine diagram. That diagram has $r_i+1$ nodes and Gram
rank $r_i$. These descriptions also apply to nonsimple ideal vertices; see
\cite{BoundsArticle,BoundsSupplement}.

Put $r=n-1$. Write the $h$-polynomial of the figure as
\begin{equation}\label{fv:hpolynomial}
 h_v(t)=\sum_{j=0}^r h_{v,j}t^j,
 \qquad
 h_v(t)=\prod_i(1+t+\cdots+t^{r_i})
 \quad\hbox{at an ideal vertex}.
\end{equation}
At a finite vertex $h_{v,j}=1$ for every $0\le j\le r$.
Define
\begin{equation}\label{fv:allowance-general}
 V_v=\sum_{j=0}^r h_{v,j},\qquad
 G_v=\sum_{j=0}^r h_{v,j}\min(j,r-j),\qquad
 b_n(v)=\frac52G_v-\frac r2V_v.
\end{equation}
The quantity $V_v$ is the number of vertices of the vertex figure.
It equals $n$ at a finite vertex and $\prod_i(r_i+1)$ at an ideal
vertex. Every vertex of the figure has degree $r$, so its number of
edges is $E_v=rV_v/2$. These edges correspond to corners at $v$.

For a finite vertex, the sum in \eqref{fv:allowance-general} gives
\begin{equation}\label{fv:finite-allowance}
 \beta_n=b_n(v)=
 \begin{cases}
 n(n-6)/8,&n\text{ even},\\
 (n-1)(n-5)/8,&n\text{ odd}.
 \end{cases}
\end{equation}
At an ideal vertex we write
\begin{equation}\label{fv:ideal-allowance}
 K_n(v)=\frac{5\sum_jh_{v,j}\min(j,n-1-j)-(n-1)V_v}{2}.
\end{equation}
Thus $K_n(v)=b_n(v)$; we use $K_n(v)$ for the bound at an ideal
vertex and $\beta_n$ for the bound at a finite vertex.

We use the following incidence inequality.

\begin{assumption}[Incidence excess]\label{fv:incidence-excess}
Let $P\subset\HH^n$ be a finite-volume Coxeter polytope, with the
vertex figures described above. The incidence theorem of
\cite{BoundsSupplement} gives
\begin{equation}\label{fv:incidence-excess-eq}
 D_5(P)\ge\sum_v b_n(v)+\Delta_n,
 \qquad \Delta_n\ge5(n-1)>0.
\end{equation}
\end{assumption}

The proof of this inequality, by counting incidences according to
height, is given in the supplement. To obtain a contradiction it
suffices to prove $D_5(P)\le\sum_vb_n(v)$; the positive term
$\Delta_n$ already gives the required strict inequality.

\subsection{The incidence linear program and its dual}
\label{fv:incidence-primal}

First fix rational functions on corners, collected in a vector
$\Psi(\kappa)\in\Q^d$. These functions will be specified in
Section~\ref{fv:features}. Their values at a corner are the same in
the face sum and in the vertex sum. A vertex type $V$ consists of a
vertex diagram, or of an ideal-link diagram, with the markings needed
to evaluate these functions. A small-face type $Q$ consists of a
triangle or quadrilateral with compatible markings at its corners.
Set
\begin{equation}\label{fv:feature-sums}
 A(V)=\sum_{\kappa\text{ at }V}\Psi(\kappa),\qquad
 A(Q)=\sum_{\kappa\text{ on }Q}\Psi(\kappa),\qquad
 d(Q)=5-s(Q).
\end{equation}
The families of types must include every vertex and every triangle
or quadrilateral of a polytope under consideration. We may also
include formal types that are not realized by a polytope.

Introduce nonnegative variables $x_V$ for vertex types, $y_Q$ for
small-face types, and $z_\kappa$ for allowed individual corners. The
last variables account for corners belonging to faces with at least
five sides. The incidence feasibility problem is
\begin{align}
 \sum_VA(V)x_V-\sum_QA(Q)y_Q-
       \sum_\kappa\Psi(\kappa)z_\kappa&=0,
       \label{fv:primal-balance}\\
 \sum_Qd(Q)y_Q-\sum_Vb_n(V)x_V&=1.
       \label{fv:primal-normalization}
\end{align}
For an infinite family of types we require only finitely many
variables to be nonzero. This suffices because a polytope has finitely
many incidences.

\begin{lemma}\label{fv:actual-lift}
A finite-volume Coxeter polytope gives a feasible point of
\eqref{fv:primal-balance}--\eqref{fv:primal-normalization} if every
vertex, every triangle or quadrilateral, and every corner has one of
the specified descriptions.
\end{lemma}

\begin{proof}
Let $f_3$ and $f_4$ be the numbers of triangular and quadrilateral
two-faces, and put
\[
 L=2f_3+f_4-\sum_vb_n(v).
\]
Since the other two-faces contribute nonpositively to $D_5(P)$,
Input~\ref{fv:incidence-excess} gives
$L\ge\Delta_n>0$. Divide each actual vertex count, small-face count
and residual-corner count by $L$. Counting the same corners by vertices
and by faces gives \eqref{fv:primal-balance}. The definition of $L$
gives \eqref{fv:primal-normalization}.
\end{proof}

The variables $z_\kappa$ cannot be omitted: corners on faces with at
least five sides still occur in the incidence equalities, although
those faces have nonpositive contributions to $D_5(P)$.

\begin{proposition}[Dual weight inequalities]\label{fv:charge-certificate}
Suppose a coefficient vector $\theta\in\R^d$ satisfies
\begin{align}
 b_n(V)-\theta\cdot A(V)&\ge0 &&\text{for every allowed vertex type},
    \label{fv:vertex-certificate}\\
 \theta\cdot A(Q)-d(Q)&\ge0 &&\text{for every allowed small face},
    \label{fv:face-certificate}\\
 \theta\cdot\Psi(\kappa)&\ge0 &&\text{for every allowed corner}.
    \label{fv:corner-certificate}
\end{align}
Then the incidence feasibility problem is infeasible.
\end{proposition}

\begin{proof}
Use Theorem~\ref{thm:homogeneous-certificate}, reversing the sign of
\eqref{fv:primal-balance}. The vertex, small-face and residual-corner
columns $(T,d)$ are respectively $(-A(V),-b_n(V))$,
$(A(Q),d(Q))$ and $(\Psi(\kappa),0)$; take $\lambda=\theta$.
The three displayed conditions are exactly its certificate inequalities.
\end{proof}

Equivalently, the weight function
$w(\kappa)=\theta\cdot\Psi(\kappa)$ has sum at least two on each
triangle, at least one on each quadrilateral, and at most $b_n(v)$
at each vertex. On a face with at least five sides its nonnegative
sum is at least $5-s(F)$. Summing over faces and vertices gives
$D_5(P)\le\sum_vb_n(v)$, contrary to
Input~\ref{fv:incidence-excess}. This is weak duality, applied to a
finite sum of incidences; no assertion of strong duality for an
infinite-dimensional problem is needed.

\subsection{The functions defining the corner weights}
\label{fv:features}

We now specify the functions used to construct the weights. Write
$\phi(\kappa)\in\Q^{26}$ for their values at a corner. Each
coordinate of $\phi$ is fixed; its coefficient in the weight function
is a variable of the linear program.

At a finite vertex, the spherical diagram is a disjoint union of
connected finite diagrams. Its two marks are the two nodes determining
the corner. At an ideal vertex a corner corresponds to an edge of
the product \eqref{fv:link}: one simplex factor varies along an edge,
and a vertex is selected in every other factor. The varying affine
diagram has two marked nodes. Each other factor contributes the finite
diagram obtained by deleting the node selected by its vertex. These
finite diagrams are included in the corner description.

We number $A_r$ along its path by $0,\ldots,r-1$. We number $D_r$
along the spine $0,\ldots,r-3$, with leaves $r-2,r-1$ attached to
$r-3$. The remaining finite and affine diagrams use the fixed
numberings of the rational tables in \cite{BoundsSupplement}.
Fix one numbering of each vertex diagram and use it at all its
corners. Under a diagram isomorphism, transport this numbering
together with the markings. The numbering is not chosen separately
at each corner.

Set $\phi_0=1$. The functions $\phi_1,\phi_2,\phi_3$ are respectively
the prescribed ordinary pair weight for two marks in the same finite
component, the affine pair weight in the varying ideal factor, and
the ordinary weight for marks in different finite component copies;
each is zero outside its stated case. Their fixed rational tables,
including exceptional small ranks, are the intrinsic corner tables of
\cite{BoundsSupplement}. We use their final edition: in particular,
the ordinary $B_{13}$ endpoint entry is zero.

The next fifteen coordinate functions are indicators. Each is one when the two
marks lie in a copy of the displayed component at the displayed
unordered pair, and when the final column holds. Otherwise it is zero.
A component in the last column is a connected component of the
corner diagram disjoint from the marked component. Put
$M_{15}=m_{A_{15}}+m_{D_{15}}$, the number of such finite components.
\begin{center}\small
\begin{tabular}{cllp{.30\textwidth}}
\toprule
$j$&Marked type&Marked pair&Other components\\\midrule
4&$A_{13}$&$\{0,12\}$&An $A_1$ occurs.\\
5&$A_{14}$&$\{0,13\}$&An $A_1$ occurs.\\
6,7&$A_{15}$&$\{0,13\},\{1,14\}$&An $A_1$ occurs.\\
8,9&$D_{15}$&$\{1,13\},\{1,14\}$&An $A_1$ occurs.\\
10,11&$A_{16}$&$\{1,12\},\{3,14\}$&An $A_1$ occurs.\\
12,13&$D_{16}$&$\{0,14\},\{0,15\}$&An $A_1$ occurs.\\
14,15&$D_{14}$&$\{0,12\},\{0,13\}$&An $A_1,A_2$ or $A_3$ occurs.\\
16,17&$D_{13}$&$\{0,11\},\{0,12\}$&$M_{15}\ge3$.\\
18&$D_{28}$&$\{0,12\}$&$M_{15}\ge2$.\\
\bottomrule
\end{tabular}
\end{center}
In a row with two indices, the two marked pairs correspond
in the displayed order.
The final seven indicators concern different finite component copies:
\begin{center}\small
\begin{tabular}{cp{.48\textwidth}p{.31\textwidth}}
\toprule
$j$&Marked components and allowed nodes&Other components\\\midrule
19&$A_4:\{0,3\}$ and $A_{12}:\{1,10\}$&An $A_1$ occurs.\\
20&$A_2:\{0,1\}$ and $A_4:\{0,3\}$&An $A_1$ and an $A_{10}$ occur.\\
21&$A_2:\{0,1\}$ and $D_5:\{0\}$&An $A_1$ and an $A_9$ occur.\\
22&$A_1:\{0\}$ and $D_5:\{0\}$&An $A_{11}$ occurs, or both an $A_9$ and an $A_2$ occur.\\
23&$A_1:\{0\}$ and $D_5:\{0\}$&An $A_9$ and an $A_2$ occur.\\
24&$A_{12}:\{0,11\}$ and $D_{16}:\{0\}$&$M_{15}\ge2$.\\
25&Two distinct $A_2$ copies, with either node in each&An $A_{10}$ and at least two $A_1$ copies occur.\\
\bottomrule
\end{tabular}
\end{center}
These functions are indicators of unordered marked pairs. Functions
22 and 23 may both be one; function 25 is zero when both marks belong
to the same $A_2$ component. Every coordinate of $\phi$ lies in
$[0,1]$, but a linear combination of them may have either sign.

There is one additional ideal-corner indicator, denoted $\chi_B$.
Write $\Gamma[m]$ for an affine diagram with $m$ nodes and rank $m-1$.
The indicator is one precisely when the varying factor is $B[13]$,
the marked pair is $\{10,11\}$ or $\{10,12\}$, and the complete
residual finite diagram contains at least three components of types
$A_{15}$ or $D_{15}$. It is zero otherwise. Thus $B[13]$ denotes
standard type $\widetilde B_{12}$, not a finite rank-13 component.
Form the residual finite diagram by deleting one selected node from
each of the other affine factors. The two marks in the varying factor
are kept; deleting them first would define a different indicator.

Let $\epsilon_n=[n=58]$. Use free vectors
$a^{58},a^{\rm tail}\in\R^{26}$ and parameters $\lambda,\rho\ge0$.
The weight function to be determined is
\begin{equation}\label{fv:piecewise-charge}
 w(\kappa)=
 \bigl(\epsilon_na^{58}+(1-\epsilon_n)a^{\rm tail}\bigr)
       \cdot\phi(\kappa)+\epsilon_n\lambda\chi_B(\kappa).
\end{equation}
The 53 functions used in the incidence equalities are
\begin{equation}\label{fv:conserved-features}
 \Psi(\kappa)=
 (\epsilon_n\phi(\kappa),(1-\epsilon_n)\phi(\kappa),
                     \epsilon_n\chi_B(\kappa)).
\end{equation}
The auxiliary variable $\rho$ gives a common nonnegative lower
bound in the vertex inequalities; there are 54 variables common to
all constraints. The normalization at an ideal vertex is specified
in Section~\ref{fv:ideal-block}. The coefficient vectors are
unrestricted in sign. We require the resulting weight to be
nonnegative at every corner in the specified classes, rather than
require every coefficient to be nonnegative.

The finite tables used for the dimension-58 calculation are defined
only through rank 58. Their validity in higher dimensions is not
assumed. The face theorems for arbitrary rank and the vertex estimates
below justify the weight function~\eqref{fv:piecewise-charge} for
$n\ge59$.

\subsection{Local face inequalities}
\label{fv:face-inputs}

The local face theorems restrict the marked diagrams and supporting
hyperplanes that can occur at a triangle or quadrilateral. They give
the following sufficient linear conditions on the weight coefficients.

\begin{assumption}[Local face theorems]\label{fv:face-assumption}
The corner, triangle and quadrilateral theorems of
\cite{BoundsSupplement}, together with the implications between
coefficient conditions established in \cite{LPRepository}, give a
linear system $\mathcal L$.
It includes the corner and triangle inequalities for both coefficient
vectors; the conditions involving $\chi_B$ and the finite $D_{13}$ conditions;
the conditions for fork diagrams and the additional rank-16 conditions;
and all the quadrilateral inequalities, including the $A_{13}$
inequality for each coefficient vector. If the coefficients and
auxiliary variables satisfy this system, then
\eqref{fv:face-certificate}--\eqref{fv:corner-certificate} hold at
every corner, triangle and quadrilateral of a polytope in dimension
$n\ge58$. Whenever one of these implications uses equalities
implied by the dimension-58 vertex inequalities, the conditions of
Section~\ref{fv:finite-network} are also required. Those equalities
are consequences of that linear system, not geometric assumptions.
\end{assumption}

The unknowns in $\mathcal L$ are coefficients of a weight function,
not quantities assigned by the polytope. Thus its inequalities are
sufficient conditions for the desired weights; they need not hold for
every weight function satisfying the face and vertex bounds. For
example, the $A_{13}$ face theorem is applied under
\begin{equation}\label{fv:A13-rows}
 a^{58}_0+a^{58}_4\ge1,
 \qquad a^{\rm tail}_0+a^{\rm tail}_4\ge1.
\end{equation}
We include both inequalities in $\mathcal L$.

The hypotheses record the identities of the facets of $P$, the order
of the marks, their embeddings in the numbered diagrams, the selected
node deletions at ideal links, and all the facets used to justify an
auxiliary construction. An inverse-Gram calculation does not establish
these geometric hypotheses. Inequalities from different alternatives
are combined only when they apply simultaneously, with each corner
counted with its stated multiplicity.

For example, for either $z=a^{58}$ or $a^{\rm tail}$, the system
contains
\begin{align}
 4z_0+z_3/16+z_{22}+z_{23}&\ge1,
       \label{fv:example-face-P}\\
 4z_0+z_3/2+z_{21}&\ge1.
       \label{fv:example-face-M}
\end{align}
Each inequality is used only for the marked face configurations in
its local theorem. The occurrence of a diagram with the same name
elsewhere on a face does not suffice. The electronic matrix lists
every inequality of $\mathcal L$ and identifies its source theorem.

\subsection{Finite vertices: a network formulation}
\label{fv:finite-network}

At a finite vertex the marks range over all unordered pairs of
nodes in a spherical forest. A formal forest has nonnegative integer
component multiplicities $m_T$ satisfying
\begin{equation}\label{fv:forest-rank}
                  \sum_T r_Tm_T=58,
\end{equation}
where $r_T$ is the rank of type $T$. We use the finite collection of
weight types of the rational catalogue. A symbolic dihedral family
can be represented by one weight type only when every corner
function has the same table throughout the stated parameter range and that
range is covered by the local theorem. A finite sample of labels is
not used in place of an unbounded family.

Let
\[
 \Phi(m)=\sum_{\kappa\text{ at the formal vertex}}\phi(\kappa).
\]
Every actual finite vertex supplies an allowed formal forest. The
converse is not required. We seek the family of inequalities
\begin{equation}\label{fv:finite-budget-endpoint}
                   a^{58}\cdot\Phi(m)+\rho\le377
                  \quad\hbox{for every }m\text{ satisfying
                  \eqref{fv:forest-rank}}.
\end{equation}
The right side is the vertex bound
$\beta_{58}=58(58-6)/8$ obtained from the incidence inequality.

We express this family of inequalities by paths in a finite directed
acyclic graph.
For one component $T$, let $S_T$ be the sum of its ordinary
same-component entries; for different component types $T,U$, let
$B_{TU}$ be the sum of their ordinary cross entries. Using $u_j$
for the $j$th coordinate vector of $\Q^{26}$ gives
\begin{equation}\label{fv:first-aggregates}
 \Phi_0=\binom{58}{2},\qquad
 \Phi_1=\sum_Tm_TS_T,\qquad
 \Phi_2=0,\qquad
 \Phi_3=\sum_{T<U}m_Tm_UB_{TU}.
\end{equation}
Functions 4--18 sum to the component multiplicity times the indicator
of the stated condition on the other components. The remaining coordinates follow by counting the unordered
mark pairs in the table: for example,
$\Phi_{19}=4m_{A_4}m_{A_{12}}[m_{A_1}>0]$ and
$\Phi_{25}=4\binom{m_{A_2}}2[m_{A_{10}}>0,m_{A_1}\ge2]$.
Brackets denote Boolean indicators; the factors four count node choices,
not orientations.

Fix the multiplicities $h_H=m_H$ of the types in
\begin{equation}\label{fv:high-set}
 \mathcal H=\{A_9,\ldots,A_{17},B_9,\ldots,B_{12},
                    D_9,\ldots,D_{16},D_{28}\}.
\end{equation}
Also fix $a=m_{A_1}$ and $b=m_{A_2}$, and denote the tuple of these
fixed multiplicities by $I$. Let
\[
 M_h=58-\sum_Hr_Hh_H,\qquad R_I=M_h-a-2b.
\]
Every nonzero ordinary cross term has exactly one component in
$\mathcal H$. Once $I$ is fixed, the remaining contribution to
$\Phi(m)$ is additive in the other component multiplicities, apart
from one indicator for the presence of $A_3$. More precisely, let
$\mathcal R$ be the remaining types after removing
$\mathcal H,A_1,A_2$, and put
$\mathcal S=\mathcal R\setminus\{A_3\}$. Define
\[
 p_T^h=S_Tu_1+\Bigl(\sum_Hh_HB_{TH}\Bigr)u_3.
\]
For remaining types other than $A_4,D_5$ use $p_{I,T}=p_T^h$.
For these two types, the vector contributions are
\begin{align}
 p_{I,A_4}={}&p_{A_4}^h+4h_{A_{12}}[a>0]u_{19}
                    +4b[a>0,h_{A_{10}}>0]u_{20},\notag\\
 p_{I,D_5}={}&p_{D_5}^h+2b[a>0,h_{A_9}>0]u_{21}\notag\\
 &+a[h_{A_{11}}>0\text{ or }(h_{A_9}>0,b>0)]u_{22}
                    +a[h_{A_9}>0,b>0]u_{23}.
                    \label{fv:residual-costs}
\end{align}
Let $K_I$ contain $\binom{58}{2}u_0$ and all terms of the corner
sum involving only the fixed counts $h,a,b$, including the fixed
same-component indicators. Put
\[
              \gamma_I=h_{D_{14}}[a=b=0](u_{14}+u_{15}).
\]
Then, with $k=m_{A_3}$,
\begin{equation}\label{fv:vector-cost-identity}
 \Phi(m)=K_I+kp_{I,A_3}+[k>0]\gamma_I
                         +\sum_{T\in\mathcal S}m_Tp_{I,T}.
\end{equation}
By the definitions of the functions, the only indicator not fixed by
$I$ is the occurrence of $A_3$ in the $D_{14}$ conditions. The
remaining cross terms are linear in the $A_4,D_5$ multiplicities;
all other terms are fixed by $I$. This proves the identity.

For each choice $(I,k)$ consistent with total rank 58, an initial
edge has vector weight $K_I+kp_{I,A_3}+[k>0]\gamma_I$ and ends
at a vertex labelled by the remaining rank $R_I-3k$. A subsequent
edge selects a type $T\in\mathcal S$, decreases the remaining rank
by $r_T$, and has vector weight $p_{I,T}$. Use the rank-feasible initial choices and their reachable remaining
ranks. A complete path ends at rank zero. The rank strictly decreases
along each subsequent edge, so the graph is acyclic.

\begin{lemma}[Exact network formulation]\label{fv:network-exact}
Every formal forest in the specified family is represented by a
complete path with vector weight $\Phi(m)$, and every complete path
gives component multiplicities of total rank 58. Consequently the
maximum scalar path weight is
$\max_m a^{58}\cdot\Phi(m)$ for every signed vector $a^{58}$.
\end{lemma}

\begin{proof}
The initial choice and subsequent edges select exactly the
nonnegative multiplicities in \eqref{fv:forest-rank}. Equation
\eqref{fv:vector-cost-identity} gives their vector weight. Different
orders may represent the same forest, but give the same weight.
A unit flow on an acyclic graph is a convex combination of complete
paths, so maximizing its linear weight gives the maximum path weight.
\end{proof}

For fixed $h$, two choices of $I$ may use the same graph of remaining
ranks only if their vectors $p_{I,A_4},p_{I,D_5}$ agree. The other
edge weights are already fixed by $h$. Their initial edges are still
kept separately. Equality of scalar weights at one choice of
coefficients would not justify identifying these graphs while the
coefficients vary.

The finite calculation is divided into subgraphs indexed by $S$.
Within one subgraph, let $g$ index the groups of initial choices with
equal remaining edge-weight vectors, as described above. Introduce free
potential variables $v_{S,g,\ell}$ at the reachable remaining ranks
$\ell$, and a free upper bound $U_S$ for the paths starting in $S$.
The linear constraints are
\begin{align}
 v_{S,g,0}&=0,\notag\\
 v_{S,g,\ell}-v_{S,g,\ell-r_T}
                  &\ge a^{58}\cdot p_{g,T},\notag\\
 U_S-v_{S,g(I),R_I-3k}
                  &\ge a^{58}\cdot(K_I+kp_{I,A_3}+[k>0]\gamma_I),
                  \label{fv:network-lp}\\
 U_S+\rho&\le377.\notag
\end{align}
All sets $S$ use the same coefficients $a^{58},\rho$ and include
every remaining type; the potentials are separate variables. Summing
along a path gives its required upper bound. Conversely, define each
potential to be the maximum scalar weight of a path from that
remaining rank to zero. These potentials satisfy the inequalities.
Eliminating them therefore gives exactly
\eqref{fv:finite-budget-endpoint}, including when some coefficients
are negative. All optimization variables are real. The exhaustive
construction of the finite graph is supplied with the checked data;
the lemma proves the equivalence with the displayed linear program.

\subsection{Ideal vertices and estimates for higher dimensions}
\label{fv:ideal-block}

Let an ideal link have factor ranks
$r_1,\ldots,r_k$, total rank $R=n-1$ and $V=\prod_i(r_i+1)$ vertices.
For a connected affine factor $\Gamma$ of rank $r$, let $Q_\Gamma$
be the sum of its baseline affine pair weights, given by $\phi_2$,
and put
$q_\Gamma=Q_\Gamma/(r+1)$. Define $C=\sum_iq_{\Gamma_i}$.
Denote the sum of the weights at $v$, before the $\chi_B$ term is
included, by $C(v)$. For coefficient values $(a_0,a_2)$ it satisfies
\begin{equation}\label{fv:ideal-charge}
              \frac{C(v)}{V}=\frac R2a_0+Ca_2
\end{equation}
This follows because a marked
pair in one factor occurs once for each choice of a vertex in every
other factor. Only coordinates 0 and 2 of $\phi$ can be nonzero at
ideal corners; $\chi_B$ is a separate function.

Let $J$ be the sum of independent random variables, with the $i$-th
variable uniformly distributed on $\{0,\ldots,r_i\}$, and put
$Y=2J-R$. The coefficient distribution
of $h_v$ is the distribution of $J$ multiplied by $V$. Symmetry gives
\begin{equation}\label{fv:normalized-ideal}
 \frac{K_n(v)}V=\frac{3R}{4}-\frac54\mathbb E|Y|,
 \qquad
 S_a=\frac{3R}{4}-\frac54\mathbb E|Y|-\frac R2a_0-Ca_2.
\end{equation}
The desired inequality is $S_a\ge\rho$. Since these expressions
have been divided by $V$, the corresponding difference between the
vertex bound and its weight sum is at least $V\rho$.

\paragraph{Products of total rank 57.}
For $n=58$, introduce free factor-bound variables $k_r$ for
$1\le r\le57$. For every connected affine family impose
\begin{equation}\label{fv:ideal-factor-links}
                 k_r\ge\frac r2a^{58}_0+q_\Gamma a^{58}_2.
\end{equation}
The family ranges are affine $A$ at ranks 1--57,
$C$ at 2--57, $B$ at 3--57, $D$ at 4--57, and the exceptional
$E_6,E_7,E_8,F_4,G_2$. There are 227 individual-family rows.
They all remain, even when $a^{58}_2$ is negative.

The ideal-vertex calculation uses 354 selected symmetric arrays. For
each array the data specify a preceding array and the rank of the
factor added to obtain it. For every other factor addition, the data
specify a selected array whose central averages are no greater than
the resulting averages. One selected index must work for all
coordinates of this comparison. These indices are fixed data; the
array averages below are variables of the linear program. Adding a
factor acts on an array by uniform convolution followed by subtraction
of its factor bound:
\[
 f'_j=\frac1{r+1}\sum_{\ell=0}^r f_{j+\ell}-k_r.
\]
The initial array is $g_j=(5\min(j,57-j)-57)/2$; after the factors
have total rank 57 it is a scalar, equal to the normalized vertex
bound minus the sum of the factor bounds. After accumulated rank
$R$, put $L=57-R$ and let $B_{R,i,a}$ be the average of array $i$
on the central interval $a,\ldots,L-a$, for
$0\le a\le\lfloor L/2\rfloor$. Let $I_R$ be the selected indices
at rank $R$. If a rank-$r$ factor is added to an array
$j\in I_{R-r}$, put $A=L-2a+1$ and
\begin{align}
 \omega_{L,a,r,\ell}&=\frac{A+r-2\ell}{A(r+1)},
       &&0\le\ell<\min(A,r+1),\notag\\
 T_{R,i,a}^{r,j}&=B_{R,i,a}+k_r-
       \sum_{\ell=0}^{\min(A,r+1)-1}
           \omega_{L,a,r,\ell}B_{R-r,j,a+\ell}.
       \label{fv:ideal-average-transform}
\end{align}
These positive coefficients sum to one and express the central average
of the convolution as an average of nested central intervals. The rows
are:
\begin{enumerate}[label=\textup{(\roman*)},nosep]
\item the 29 initial equalities
\[
 B_{0,0,a}=\frac1{58-2a}\sum_{j=a}^{57-a}
                  \frac{5\min(j,57-j)-57}{2},
 \qquad0\le a\le28;
\]
\item $T_{R,i,a}^{r,j}=0$ for the specified predecessor $(r,j)$ of
 each selected $i\in I_R$;
\item $T_{R,D(R,r,j),a}^{r,j}\le0$ for every candidate, using
 its fixed selected comparison array $D(R,r,j)$;
\item $\tau\le B_{57,i,0}$ for every final selected index, and
 $\rho\le\tau$, with $\tau$ free.
\end{enumerate}
The index sets and maps are included in the electronic ideal-vertex
data of \cite{LPRepository}. Together with the connected-family
inequalities, this system has 103,265 rows and 5,754 auxiliary
variables. Positivity of convolution preserves the comparisons at
every factor addition. Hence every product of total rank 57 has the
required difference at least $\tau\ge\rho$. The fixed choices of
comparison arrays give a sufficient system of linear inequalities;
they need not characterize every valid ideal-vertex weight function.

\paragraph{The additional ideal-corner indicator.}
Keep the preceding additive inequalities and treat $\chi_B$
separately. This indicator depends on the other factors as well as
the rank-12 factor, so its contribution cannot be included in an
unconditional change of $k_{12}$.
If $\chi_B$ occurs, the varying $B[13]$ factor has rank 12 and
the residual finite rank is 45. Three rank-15 components exhaust that
rank. The complete one-node-deletion classification leaves exactly
\begin{align}
 \mathrm{I}:&\quad B[13]\times\Gamma_1[16]
                   \times\Gamma_2[16]\times\Gamma_3[16],\notag\\
 \mathrm{II}:&\quad B[13]\times D[31]\times\Gamma[16],
 \qquad\Gamma_i,\Gamma\in\{A,B,D\}.
 \label{fv:guarded-links}
\end{align}
Let $a_\chi=V^{-1}\sum_\kappa\chi_B(\kappa)$. The additional rows are
\begin{equation}\label{fv:guarded-rows}
  \frac{57}{2}a^{58}_0+C a^{58}_2+a_\chi\lambda+\rho\le K.
\end{equation}
All thirteen rows have the following compact specification. Put
\[
 (c_A,c_B,c_D)=(5,25/4,411/64),\qquad
 (d_A,d_B,d_D)=(1,1/16,1/4).
\]
For every unordered triple $\Gamma_1,\Gamma_2,\Gamma_3\in\{A,B,D\}$
in type I, and every $\Gamma\in\{A,B,D\}$ in type II, use
\begin{equation}\label{fv:guarded-data}
\begin{array}{c|ccc}
 &C&a_\chi&K\\\hline
 \mathrm{I}&11/2+\sum_{i=1}^3c_{\Gamma_i}
       &(2/13)\prod_{i=1}^3d_{\Gamma_i}&331643/13312\\
 \mathrm{II}&753/62+c_\Gamma&(2/403)d_\Gamma&16549/806.
\end{array}
\end{equation}
Here $d_\Gamma$ is the proportion of one-node deletions giving the
required residual component; the deletion counts are $16,1,4$.
The values of $K$ are \eqref{fv:normalized-ideal} for ranks
$(12,15,15,15)$ and $(12,30,15)$. Thus the specification retains
deletion multiplicities, not just the affine factor types.

\paragraph{The remaining rank-58 products in dimension 59.}
Append all 63 connected-family links of ranks at most 16 for the
coefficients $a^{\rm tail}$,
\begin{equation}\label{fv:tail-factor-links}
             k_r\ge\frac r2a^{\rm tail}_0
                                      +q_\Gamma a^{\rm tail}_2,
\end{equation}
using the \emph{same} factor bounds $k_r$ as above. Put $k_0=0$ and introduce
$\mu\ge0$. The fifteen descent rows are
\begin{equation}\label{fv:descent}
 k_r-k_{r-1}+\mu\le h_r,\qquad 1\le r\le15,
\end{equation}
where
\[
 (h_1,\ldots,h_{15})=
 \left(\frac34,\frac{47}{66},\frac{89}{132},\frac7{11},
 \frac{79}{132},\frac{37}{66},\frac{23}{44},
 \frac{2213}{4374},\frac{59}{132},\frac9{22},
 \frac{67}{156},\frac{71}{168},\frac5{12},
 \frac{79}{192},\frac{83}{204}\right).
\]
The rank-reduction theorem applies to a rank-58 product whose largest
factor rank is between 10 and 16. It provides a factor whose rank
can be reduced by one, with the difference between the normalized
vertex bound and the sum of factor bounds increasing by at least
$\mu$ when the factor is restored. The reduced product has rank 57
and satisfies the preceding inequalities. Thus the dimension-59
vertex inequality holds for this class. The same variables $k_r$
must be used in both systems for the comparison to be valid.

\paragraph{All other ideal links in dimensions at least 59.}
For every other dimension-59 ideal link, and for every ideal link
in dimension at least 60, the analytic vertex theorem gives
$S_{(0,1)}\ge\varepsilon=1/5307$. For coefficients
$a^{\rm tail}$ the following four affine inequalities suffice:
\begin{align}
 a^{\rm tail}_0&\le0,&
 a^{\rm tail}_0+a^{\rm tail}_2&\le1,\notag\\
 \rho+29a^{\rm tail}_0&\le\varepsilon,&
 \rho+29a^{\rm tail}_0+29a^{\rm tail}_2&\le29+\varepsilon.
 \label{fv:analytic-ideal-rows}
\end{align}
To see the reduction, write $t=1/n$,
$u_0=a^{\rm tail}_0$, $u_1=a^{\rm tail}_0+a^{\rm tail}_2-1$ and
\[
 P_i(t)=(\varepsilon-\rho)t-\frac{1-t}{2}u_i.
\]
Each $P_i$ is affine, and \eqref{fv:analytic-ideal-rows} is exactly
nonnegativity at the two endpoints of $[0,1/59]$. The analytic input
also gives $0\le C\le R/2$. If $w=2C/R$, then
\[
 t(S_a-\rho)=t(S_{(0,1)}-\varepsilon)
                         +(1-w)P_0(t)+wP_1(t)\ge0.
\]
Thus the four inequalities prove the vertex bound throughout the
stated range. This conclusion uses the analytic estimate for
$S_{(0,1)}$ in every dimension of that range.

\paragraph{Finite vertices in dimensions at least 59.}
Let $q\in\Q^{26}$ be the vector displayed in
\eqref{fv:successful-law} below. Introduce $0\le\xi\le1$ and impose
\begin{equation}\label{fv:finite-tail-rows}
 a^{\rm tail}_j\le q_j+2\xi[\,j\in\{16,17\}\,]
       \quad(j\ne2),\qquad \rho\le\frac14.
\end{equation}
The analytic finite-vertex estimate asserts that, for every
formal finite forest of every rank $n\ge59$,
\[
 \beta_n-\bigl(q+2\xi(u_{16}+u_{17})\bigr)\cdot\Phi_n
                                                       \ge\frac14.
\]
Since the functions at a finite vertex are nonnegative and function
2 is zero, \eqref{fv:finite-tail-rows} proves its vertex bound. The
analytic estimate permits every other component and every rank
$n\ge59$; it is separate from the rank-58 finite calculation. The
coordinate inequalities are sufficient conditions on
$a^{\rm tail}$, not a characterization of all valid finite-vertex
weight functions.

\subsection{The optimal bound within the chosen weight family}
\label{fv:377}

Fix a finite family $\mathcal T$ of formal forests and a family
$\mathcal C$ of weight coefficients satisfying the local face and
ideal-vertex conditions. Define the maximum finite-vertex weight sum
and minimize it over $\mathcal C$:
\begin{equation}\label{fv:budget-template}
 U_{\mathcal T}(a)=\max_{F\in\mathcal T}a\cdot\Phi(F),
 \qquad \min_{a\in\mathcal C}U_{\mathcal T}(a).
\end{equation}
Here $a$ is the coefficient vector used at a finite vertex, and
$\Phi(F)$ is the sum of its corner-function vectors. The set
$\mathcal C$ is the projection onto $a$ of the joint system, including
all auxiliary-variable conditions. The incidence inequality supplies
an independent upper bound for comparison. To prove the dimension
bound, one needs a feasible choice whose maximum finite-vertex
weight sum does not exceed that upper bound, with the face and
ideal-vertex inequalities satisfied by the same coefficients.

In this application $\mathcal T$ is the exhaustive family of allowed
rank-58 forests and $\beta_{58}=377$. We give both a feasible weight
function and a dual lower bound for this coefficient family. The
rank-58 calculation is used only at that rank.

The following rational point uses the coordinate order specified in
Section~\ref{fv:features}:
\begin{align}
 (q_0,\ldots,q_{12})={}&
 (0,1,1,1,1,\tfrac12,\tfrac14,\tfrac14,\tfrac14,\tfrac14,
       \tfrac58,\tfrac58,\tfrac12),\notag\\
 (q_{13},\ldots,q_{25})={}&
 (\tfrac12,\tfrac12,\tfrac12,-\tfrac12,-\tfrac12,\tfrac12,
       \tfrac18,\tfrac14,\tfrac12,\tfrac38,\tfrac58,\tfrac12,
       \tfrac14),\notag\\
 a^{58}={}&q,\qquad
 a^{\rm tail}=q+\tfrac12(u_{16}+u_{17}),\notag\\
 \lambda={}&\tfrac12,\qquad\rho=0.
 \label{fv:successful-law}
\end{align}
The two $D_{13}$ coefficients are negative in dimension 58 and zero
in higher dimensions. Nevertheless, their sum with the other terms
in~\eqref{fv:piecewise-charge} is nonnegative at every allowed
marked corner.

The auxiliary variables have exact values as well. For the
higher-dimensional finite-vertex and rank-reduction inequalities take
$\xi=1/4$ and $\mu=13/2187$. Take
$k_r=\max_{\operatorname{rank}\Gamma=r}q_\Gamma$ over all the
connected affine families, and assign the $B_{R,i,a}$ using the
specified ideal arrays and their defining equalities, with $\tau=0$.
The rational verification checks all these variables and every
finite-network potential at the same coefficient values.

For example, the face rows
\eqref{fv:example-face-P}--\eqref{fv:example-face-M} take values
$17/16$ and $1$ for both coefficient vectors. The ideal-indicator inequalities have slack
at least $939/13312$ in type I and $50623/25792$ in type II. These
examples illustrate the assignment; rational verification checks
every local and network inequality at the same coefficients.

There is a useful optimization statement in addition to feasibility.
Let $\mathcal A$ be the projection onto
$(a^{58},a^{\rm tail},\lambda)$ of the complete local coefficient
system, including \eqref{fv:A13-rows}, with $\rho=0$. Replace the
last network row by $U_S\le U$, for one free common variable $U$,
and minimize $U$. Eliminating the potentials gives
\begin{equation}\label{fv:minimax}
 \min_{(a^{58},a^{\rm tail},\lambda)\in\mathcal A}
                 \max_F a^{58}\cdot\Phi(F),
\end{equation}
where $F$ ranges over the exhaustive family of allowed rank-58
forests. In this optimization, 377 is no longer imposed as an upper
bound.

\begin{proposition}[Exact restricted minimax]\label{fv:exact-minimax}
The optimum in \eqref{fv:minimax} is 377.
\end{proposition}

\begin{proof}
The feasible assignment \eqref{fv:successful-law}, including
its auxiliary local variables and all network potentials, gives an upper
bound of 377. For the converse, write $a_0,a_1$ for coordinates 0 and
1 of $a^{58}$. The corner and triangle inequalities give
$a_0\ge0$ and $3a_0+2a_1\ge2$.
The formal forest
\[
                    F_*=2D_{17}+D_{24}
\]
has total rank 58. The ordinary table gives
$S_{D_{17}}=227/2$ and
$S_{D_{24}}=8\cdot24-32-1-4-3-2=150$. All cross-component weights and indicator functions vanish on this
forest, so its vector of corner sums is
$\Phi(F_*)=1653u_0+377u_1$. Thus one path bound gives
$U\ge1653a_0+377a_1$. The exact nonnegative combination
\begin{equation}\label{fv:minimax-identity}
 U-377=(U-1653a_0-377a_1)
 +\frac{377}{2}(3a_0+2a_1-2)+\frac{2175}{2}a_0
\end{equation}
proves the lower bound: the first term is the sum of this forest's
network rows and its link to $U$, and the other terms are multiples
of the two local rows.
\end{proof}

The incidence inequality gives the bound 377 independently of this
optimization. Feasibility and \eqref{fv:minimax-identity} prove that
377 is also the optimum for the chosen weight family and formal
forests. The forest $F_*$ need not occur at a polytope vertex. A larger
family of weight functions may have a smaller optimum.

\begin{remark}[Choice of the weight family]\label{fv:vinberg-weight-obstruction}
Vinberg~\cite[\S6.12, p.~91]{Vinberg1984} gives a related $E_7$
obstruction for a restricted family of nonnegative corner weights in
the compact case. His requirement that weights vanish when the marks
lie in different connected components is absent from the family in
Section~\ref{fv:features}, which allows cross-component and contextual
contributions. Such obstructions concern the specified diagram and
weight classes and do not preclude stronger certificates using
additional geometric compatibility.
\end{remark}

In the system with upper bound 377,
$\max_Fa^{58}\cdot\Phi(F)+\rho\le377$ and $\rho\ge0$.
The same inequalities force $a_0=0$, $a_1=1$ and $\rho=0$.
Consequently a requirement $\rho>0$ would exclude this valid
solution. The strict incidence inequality already gives the
contradiction when $\rho=0$.

\subsection{The finite-volume dimension bound}
\label{fv:all-dimensions}

The estimates used in each dimension are summarized below. In every
case the same weight function is used for the face and vertex sums.
\begin{center}\small
\begin{tabular}{p{.29\textwidth}p{.61\textwidth}}
\toprule
Class&Estimate used\\\midrule
Finite vertices, $n=58$&The finite directed graph and all its
 potential inequalities.\\
Finite vertices, $n\ge59$&The analytic finite-vertex estimate and
 coordinate inequalities \eqref{fv:finite-tail-rows}.\\
Ideal vertices, $n=58$, no corner with $\chi_B=1$&The system for
 rank-57 products and the connected-family inequalities.\\
Ideal vertices, $n=58$, some corner with $\chi_B=1$&The exhaustive
 classification of node deletions \eqref{fv:guarded-links} and all thirteen
 rows \eqref{fv:guarded-rows}.\\
Ideal vertices, $n=59$, largest factor rank 10--16&The same factor
 bounds and the rank-reduction inequalities.\\
Other ideal vertices at 59, and all $n\ge60$&The analytic ideal
 estimate and the four inequalities \eqref{fv:analytic-ideal-rows}.\\
All corners and small faces, $n\ge58$&The exhaustive local
 geometric classification and theorems in Input~\ref{fv:face-assumption}.\\
\bottomrule
\end{tabular}
\end{center}

\begin{theorem}[Finite-volume bound]\label{fv:recovery}
Assume the local geometric estimates, the exhaustive diagram
classifications and the incidence inequality stated above. The
verified rational solution excludes finite-volume hyperbolic Coxeter
polytopes in every integer dimension $n\ge58$.
\end{theorem}

\begin{proof}
Suppose such a polytope exists. The last row of the table gives
nonnegative corner weights whose sums are at least two on each
triangle and at least one on each quadrilateral. The preceding rows give
\[
 \sum_{\kappa\text{ at }v}w(\kappa)\le b_n(v)
\]
at every finite or ideal vertex. The dual inequalities contradict
the feasible point supplied by Lemma~\ref{fv:actual-lift}.
\end{proof}

\section{Linear programs for squared projection norms}
\label{sec:energy-framework}

The two compact cases use the same specialization of
Theorem~\ref{thm:homogeneous-certificate}. For each configuration
class $b$, let $e_b$ be the vector of squared projection norms defined
in~\eqref{eq:component-energies}, grouping orthogonal components
when convenient. Write all the local inequalities valid on this
class as
\begin{equation}\label{eq:branch-energy}
                  A_be_b\ge c_b,\qquad e_b\ge0.
\end{equation}
The matrices $A_b$ and vectors $c_b$ are rational and nonnegative,
and every row of $A_b$ has a positive entry. In our applications
$A_b$ has entries $0$ and $1$: a row specifies the components in the
sum it estimates. Allowing general nonnegative coefficients gives
the same argument.

Suppose more generally that the geometric upper bound is
\begin{equation}\label{eq:general-energy-budget}
              h_b^Te_b<B_b,\qquad h_b>0,\quad B_b>0.
\end{equation}
Here $h_b$ and $B_b$ are fixed rational data. For the compact cases
$h_b=\one$ and $B_b=1$, by~\eqref{eq:component-energies}.
Variables belonging to different configuration classes are kept
separate. Within one class, every inequality concerning a given
orthogonal component uses the same coordinate of $e_b$.

\subsection{Normalization and duality}

For each class introduce nonnegative variables $\tau_b$, $E_b$ and
$\delta_b$, and impose
\begin{equation}\label{eq:budget-deficit}
 \begin{gathered}
 A_bE_b\ge c_b\tau_b,\qquad
 \delta_b=B_b\tau_b-h_b^TE_b\ge0,\\
                         \sum_b\delta_b=1.
 \end{gathered}
\end{equation}
The scalar $\tau_b$ is a scaling coefficient, $E_b$ is the scaled
norm vector, and $\delta_b$ is its difference from the scaled upper
bound. For a class ruled out by a local geometric theorem, add
\begin{equation}\label{eq:excluded-mass}
                              \tau_b=0.
\end{equation}
This also forces $E_b=0$ and $\delta_b=0$, since $h_b>0$.
The normalization by $\sum_b\delta_b$ expresses the strict
inequality as a closed linear feasibility problem, without choosing
a fixed positive error tolerance.

\begin{theorem}[Normalization and linear programming duality]
\label{thm:budget-lift}\label{prop:general-energy-budget}
\label{thm:energy-normalization}\label{thm:energy-certificate}
\label{prop:capacity-complete}
For a finite exhaustive family of configuration classes, the
following assertions hold.
\begin{enumerate}[label=\textup{(\roman*)}]
\item System~\eqref{eq:budget-deficit}--\eqref{eq:excluded-mass} is
feasible if and only if some class not locally excluded admits a
vector $e_b$ satisfying~\eqref{eq:branch-energy}
and~\eqref{eq:general-energy-budget}. In particular, a polytope in
any of the specified dimensions determines a feasible point.
\item Put
\begin{equation}\label{eq:branch-minimum}
 m_b=\min\{h_b^Te:e\ge0,\ A_be\ge c_b\}.
\end{equation}
Then
\begin{equation}\label{eq:branch-duality}
 m_b=\max\{c_b^T\alpha:\alpha\ge0,\ A_b^T\alpha\le h_b\},
\end{equation}
and both optima are attained at rational points. The normalized
system is infeasible if and only if $m_b\ge B_b$ for every class
not locally excluded.
\item Therefore, nonnegative rational multipliers satisfying
\begin{equation}\label{eq:general-capacity-certificate}
               A_b^T\alpha_b\le h_b,\qquad
               c_b^T\alpha_b\ge B_b
\end{equation}
for every such class prove infeasibility. Conversely, if the
normalized system is infeasible, these multipliers exist.
\end{enumerate}
\end{theorem}
\begin{proof}
Given a vector with $s=B_b-h_b^Te_b>0$, take
$\tau_b=1/s$, $E_b=e_b/s$, $\delta_b=1$, and set all variables
on other classes to zero. Conversely, a normalized feasible point
has some $\delta_b>0$. Then $\tau_b>0$, and division by
$\tau_b$ gives the local inequalities and
\[
           h_b^T(E_b/\tau_b)=B_b-\delta_b/\tau_b<B_b.
\]
Every polytope belongs to a class of the exhaustive family, so the
first construction applies to its squared projection norms. The scaling
coefficients need no upper bound: a small positive value of $s$
may require a large $\tau_b$.

For the minimum, a sufficiently large multiple of $\one$ satisfies
all local inequalities. Since $h_b>0$, every nonempty finite
objective sublevel is compact, so the minimum is attained. To
establish duality, omit the subscript and consider the polyhedral
cone
\[
       K=\{(Ae-u,h^Te+t):e,u,t\ge0\}.
\]
It contains $(c,m)$ but no $(c,m-\varepsilon)$ for
$\varepsilon>0$. A nonzero supporting functional $(-\alpha,\beta)$,
nonnegative on $K$ and zero at $(c,m)$, gives
\[
 \alpha\ge0,\quad\beta\ge0,\quad A^T\alpha\le\beta h,
                  \quad\alpha^Tc=\beta m.
\]
If $\beta=0$, the positive entry in each nonnegative row of $A$
forces $\alpha=0$, a contradiction. Dividing by $\beta$ gives a
dual optimizer. The reverse inequality follows by multiplying
$Ae\ge c$ by any feasible dual vector. A vertex of the compact
primal optimal face is a rational vertex of the feasible polyhedron.
The dual optimal face is then a nonempty rational polyhedron and
has a rational point by Lemma~\ref{lem:cone-alternative}.

Finally, multiplying the scaled inequalities by the dual
multipliers gives
\[
 h_b^TE_b\ge\alpha_b^TA_bE_b
            \ge(c_b^T\alpha_b)\tau_b\ge B_b\tau_b.
\]
Hence every $\delta_b$ is zero, contradicting their normalization.
In Theorem~\ref{thm:homogeneous-certificate}, take $C_b$ to be the
local cone defined by~\eqref{eq:budget-deficit}, with zero scaling
coefficient on excluded classes, $T_b=0$ and $d_b=\delta_b$.
The displayed inequality gives its dual condition $-d_b\ge0$.
Part (ii) proves the converse in part (iii).
\end{proof}

\subsection{The dual inequalities in the compact cases}

For $h_b=\one$, $B_b=1$, the computed dual solution satisfies
\begin{equation}\label{eq:capacity-dual}
 \alpha_b\ge0,\qquad A_b^T\alpha_b\le\one,\qquad
              c_b^T\alpha_b\ge\eta,\qquad\eta\ge1.
\end{equation}
For a row supported on $S_\ell$, the coefficient of $e_{bj}$ in
the chosen combination is
$\sum_{\ell:j\in S_\ell}\alpha_{b\ell}$. The inequality
$A_b^T\alpha_b\le\one$ requires this coefficient to be at most
one for each component. Nonnegativity of the unused terms then
bounds the combination by the total squared norm.

\begin{example*}
Suppose $e_I+e_A\ge a$ and $e_I+e_B\ge b$. Adding both
inequalities with coefficient one would give coefficient two on
$e_I$. Coefficients $1/2,1/2$ give instead
\[
           e_I+\tfrac12e_A+\tfrac12e_B\ge\tfrac12(a+b).
\]
The left side is at most $e_I+e_A+e_B$, since all three variables
are nonnegative. If the right side is at least one, this contradicts
the strict upper bound on the total squared norm.
\end{example*}

For a fixed dimension-labelled family, define
\begin{equation}\label{eq:threshold-energy-value}
 \Gamma_N=\min_{\substack{b\ \text{not locally excluded}\\ n_b\ge N}}m_b,
 \qquad\min\varnothing=+\infty.
\end{equation}
The normalized LP is infeasible if and only if $\Gamma_N\ge1$.
Increasing $N$ removes classes, so $\Gamma_N$ cannot decrease.
This is the specialization of Proposition~\ref{prop:threshold}.
The equivalence concerns the specified linear inequalities; a
feasible vector need not be realized by a geometric configuration.

The computation imposes~\eqref{eq:capacity-dual} and maximizes the
sum of the local dual objective values. The resulting rational
solution is checked exactly. Only its feasibility is needed for the
proof. The compact and compact arithmetic cases use different
configuration classes and local inequalities in this same program.

\section{The compact case}
\label{c:section}

We apply the local estimates of Supplement~C of~\cite{BoundsSupplement}
to the squared norms of the orthogonal projections of $x$.
These estimates hold under the hypotheses stated below.

Fix the farthest vertex and the vector $x$ defined in
\eqref{eq:original-vector}. For a connected spherical
component $C$ of the Coxeter diagram at that vertex, let
$\langle C\rangle$ denote the span of its facet normals, and put
\[
 e_C=\norm{\operatorname{proj}_{\langle C\rangle}x}^2.
\]
The spans of the components are mutually orthogonal and together form
the tangent space at the vertex. Hence
\begin{equation}\label{c:budget}
 e_C\ge0,\qquad \sum_C e_C=\norm{x}^2<1.
\end{equation}
For disjoint collections $\mathcal S$ of components, write
$e_{\mathcal S}=\sum_{C\in\mathcal S}e_C$. All such quantities in this
section refer to the same vector $x$. An estimate said to allow
arbitrary other components imposes no further condition on the remaining
spherical components at this vertex; the opposite facets are still
facets of the same polytope.

Vinberg's theorem~\cite{Vinberg1984} excludes compact Coxeter polytopes
in dimensions at least $30$. We therefore consider each dimension
$18\le n\le29$ directly.

\subsection{Small components and exceptional wall configurations}
\label{c:small}

We use the usual notation for irreducible spherical Coxeter diagrams.
The rank-two diagrams are $I_2(m)$, where $m\ge3$ is an integer;
$A_2$ denotes $I_2(3)$. The isomorphic types $B_2$ and $D_3$ are not
listed separately. A component is \emph{ordinary} if it does not belong
to one of the exceptional configurations defined below.

Each exceptional configuration consists of components of the diagram
at the chosen vertex together with a specified opposite facet. Thus
$O_7$, for example, requires more than the occurrence of a bond labelled
seven. Diagram labels always describe pairs of facets of the polytope;
projected auxiliary vectors are not assigned such labels.

\begin{proposition}[Local estimates for ordinary components]
\label{c:ordinary-input}
For each ordinary component $C$ in Table~\ref{c:ordinary-table},
$e_C$ is strictly greater than the listed bound. These estimates, proved
in Supplement~C, allow arbitrary other finite components.
\end{proposition}

\begin{table}[htbp]
\centering
\caption{The ordinary component lower bounds. Each decimal in this paper
is an exact rational number.}\label{c:ordinary-table}
\begin{tabular}{ll@{\qquad}ll}
\toprule
Component&Lower bound&Component&Lower bound\\\midrule
$A_1$&$1/10$&$A_2$&$319/2000$\\
$I_2(m),\ m\ge4$&$1/7$&$A_3$&$1/4$\\
$B_3,H_3$&$2/7$&$A_4$&$2/7$\\
$B_4,D_4,F_4$&$1/3$&$H_4$&$11/36$\\
$A_5$&$3/8$&$D_5$&$5/14$\\
$A_6,B_5,B_6,D_6$&$3/7$&$B_7,D_7$&$1/2$\\
$B_8,D_8$&$4/7$&$E_6$&$9/20$\\
$E_7$&$5/9$&$E_8$&$2/3$\\\bottomrule
\end{tabular}
\end{table}

The table includes the improved estimate for $B_3,H_3$ from
Supplement~C. The dihedral estimate holds for every integer $m\ge4$.

We next define the exceptional configurations. Let $p$ be the unit
timelike vector representing the chosen vertex and let $u$ be the unit
normal of an incident $A_1$ facet. An opposite facet whose spatial
normal lies in this rank-one component has, in the short case, a unit
normal of the form
\begin{equation}\label{c:I-normal}
 c_I=\sqrt h\,p-\sqrt{1+h}\,u,
 \qquad 0<h<\frac12.
\end{equation}
We denote this configuration by $I$.
Configuration $J$ involves two $A_1$ components and one opposite
facet whose spatial normal lies in their direct sum and is orthogonal
to every other component. That facet is separated from one of the two incident facets
(the corresponding diagram edge is dotted), and it meets the other
at angle $\pi/3$. The square of the coefficient of $p$ in its unit
normal lies in $(1/3,3/8)$. Both $A_1$ components occur together in
this configuration.

The two rank-two exceptional configurations are
\[
 \begin{aligned}
 N_7&=(\text{incident }A_2,\ \text{isolated opposite contact }7),\\
 O_7&=(\text{incident }I_2(7),\ \text{isolated opposite contact }3).
 \end{aligned}
\]
Here the contact label records the angle between the opposite facet
and an incident facet in the indicated component. Isolation means
that the spatial part of the opposite normal is orthogonal to every
other component. In particular, the definition specifies an opposite
facet as well as the incident rank-two diagram.

\begin{proposition}[Exceptional configurations]
\label{c:defect-input}
At most one of the exceptional configurations $I,J,N_7,O_7$ occurs at the
chosen vertex. Supplement~C gives the following bounds on the sum of
$e_C$ over the indicated incident components:
\[
\begin{array}{c|c|c}
\text{configuration}&\text{incident components}&\text{lower bound}\\\hline
I&A_1&0\\
J&A_1+A_1&17/50\\
N_7&A_2&9/100\\
O_7&I_2(7)&9/100.
\end{array}
\]
The positive bounds are strict. The ordinary-component bounds are not
also imposed on components belonging to an exceptional configuration.
\end{proposition}

In sums of diagrams, each of the symbols $I,J,N_7,O_7$ stands for
its indicated incident components, together with the stated conditions
on an opposite facet. Thus the incident components form the vertex
diagram; the additional facet is not included in its rank.

The assertion that these configurations cannot occur together is part
of the local geometric theorem. It includes the cases in which two
candidate opposite facets would involve a common incident component.
It does not follow from the list of spherical diagram types alone.

Supplement~C also gives estimates involving several components.
If $T_j$ is an ordinary irreducible component of rank $j$, then
\begin{align}
 e_{I+A_1}&>1/5,&
 e_{I+2A_1}&>3/10,\label{c:I-joint-one}\\
 e_{I+T_2}&>21/100,&
 e_{I+T_3}&>2/7.\label{c:I-joint-two}
\end{align}
Each left side is the sum over the components displayed in its
subscript. If several of these inequalities are used together, the
coefficient of $e_I$ must satisfy the same dual inequality in
Theorem~\ref{thm:energy-certificate}.

If $N_7$ or $O_7$ occurs, the squared projection norm for every
other ordinary rank-two component is greater than $1/5$. For an
ordinary $A_1$ component outside the exceptional configuration, it
is greater than $1/5$ in the $N_7$ case and greater than $2/15$ in
the $O_7$ case.
These bounds concern only the additional component. Their hypotheses,
and those of the preceding inequalities, are stated in the
ordinary-component appendix of Supplement~C.

\subsection{Large components}
\label{c:large}

The $A$-component theorems give
\begin{equation}\label{c:A-thresholds}
 e_{A_r}>\frac12\quad(r\ge7),\qquad
 e_{A_r}>\frac47\quad(r\ge8),\qquad
 e_{A_r}>\frac57\quad(r\ge13).
\end{equation}
All three estimates allow arbitrary other finite components.
Table~\ref{c:A-table} gives the additional estimates at specified
ranks. Each is a bound for $e_C$ of the fixed vector $x$, including
when its proof uses auxiliary points or adjacent vertices.

\begin{table}[htbp]
\centering
\caption{Additional $A$-component bounds, allowing arbitrary other finite components.}
\label{c:A-table}
\begin{tabular}{cl}
\toprule
Rank&Lower bound\\\midrule
12&$13/20$\\
15&$413/500$\\
16&$43/50$\\
17&$6/7$\\
18&$91/100$\\
19&$563/625$\\
20&$.92274$\\\bottomrule
\end{tabular}
\end{table}

For $B_r,D_r$ components with $r\ge9$, one has
\begin{equation}\label{c:BD-nine}
                         e_{B_r},e_{D_r}>\frac9{14}.
\end{equation}
The next estimate concerns an initial sequence of coordinates in a
$B_r$ or $D_r$ component, in the coordinate convention of Supplement~C.

\begin{proposition}[An estimate for initial coordinates]
\label{c:prefix-input}
Let $6\le M\le16$. Suppose the component is $B_r$ with
$r\ge M+5$, or $D_r$ with $r\ge M+6$. Put
\[
 f_i=\sqrt i-\frac{61}{100\sqrt2},\quad
 E_M=\frac{f_5^2}{5}+
        \sum_{i=5}^{M-1}(\sqrt{i+1}-\sqrt i)^2,
 \quad H_M=\left\lfloor\frac{M-5}{2}\right\rfloor.
\]
The sum of squares of the first $M+1$ coordinates of the projection
of $x$ to this component is greater than
\begin{equation}\label{c:prefix-formula}
                  E_M-H_M/3200+3969/800000.
\end{equation}
The estimate allows arbitrary other finite components.
\end{proposition}

The first two terms give the earlier estimate on the same
coordinates. The last term improves that estimate under the same
hypotheses; it is not a contribution from another subspace.
The program uses the rational lower bounds $P_M+3969/800000$,
where
\begin{center}
\begin{tabular}{c|ccccc}
$M$&7&8&9&10&11\\\hline
$P_M$&$.7349975$&$.7683275$&$.797425$&$.823735$&$.8472225$\\[2pt]
\end{tabular}
\par\smallskip
\begin{tabular}{c|ccccc}
$M$&12&13&14&15&16\\\hline
$P_M$&$.8689525$&$.88864$&$.90715$&$.9240775$&$.9401975$
\end{tabular}
\end{center}
Supplement~C proves $P_M<E_M-H_M/3200$. The finite linear program
uses $7\le M\le16$. If both the earlier and improved estimates are
included for the same coordinates, the earlier inequality is redundant.

Three additional estimates apply to the following full vertex
diagrams. They are stated separately because the preceding rank
conditions do not cover the relevant terminal coordinates:
\begin{enumerate}[label=\textup{(\roman*)}]
\item If the full vertex diagram is $B_{16}+T_2$ or $D_{16}+T_2$ in
 dimension $18$, with $T_2$ ordinary, the squared projection norm on
 the large component is greater than $.8689525$.
\item If the full vertex diagram is $B_{17}+A_1$ in dimension $18$,
 with ordinary nonshort $A_1$, the squared projection norm on $B_{17}$ is
 greater than $.90715$.
\item If the full vertex diagram is $D_{19}+A_1$ in dimension $20$,
 with ordinary nonshort $A_1$, the extension through the terminal
 $D_4$ subdiagram and the improved initial-coordinate estimate give
 $e_{D_{19}}>.91211125$.
\end{enumerate}
These local results are proved in Supplement~C at the stated ranks.
They are used as component estimates, independently of the eventual
dimension exclusions. In particular, the estimate requiring ordinary
$T_2$ is not applied to $N_7$ or $O_7$, and the ordinary $A_1$
estimates are not applied to $I$.

\subsection{An exhaustive family of configurations}
\label{c:cover}

Call a component \emph{large} if it is of type $A_r$, $r\ge7$,
or $B_r,D_r$, $r\ge9$. Every large component has squared projection norm
greater than one half by \eqref{c:A-thresholds} and
\eqref{c:BD-nine}. Table~\ref{c:ordinary-table} shows that every
other ordinary component contributes at least its rank divided by
fourteen.

\begin{lemma}[A bound for the remaining components]\label{c:spectator}
Let $\mathcal S$ be a collection of components, none large, of total
rank $S$. Whenever an exceptional configuration occurs, include all
of its incident components in the collection. Then
\begin{equation}\label{c:spectator-rate}
                          e_{\mathcal S}\ge\frac{S-1}{14}.
\end{equation}
If no exceptional configuration occurs, or the configuration is $J$,
the stronger bound $S/14$ holds. The possible loss is $1/14$ in
case $I$ and $37/700$ in cases $N_7$ and $O_7$.
\end{lemma}
\begin{proof}
Sum the rank-over-fourteen estimates on the ordinary components.
Configuration $I$ replaces a rank-one contribution by a nonnegative
quantity, so decreases the sum by at most $1/14$. The two incident
components of $J$ have total rank two, and the bound $17/50$ exceeds
$2/14$. Each rank-two configuration $N_7$ or $O_7$ contributes at
least $9/100$; the difference from $2/14$ is $37/700<1/14$.
Proposition~\ref{c:defect-input} permits at most one exceptional
configuration.
\end{proof}

Fix $18\le n\le29$. We construct an exhaustive family of
configuration classes as follows. A class may specify only component
ranks and norm bounds when the individual diagram types are not needed.

\begin{enumerate}[label=\textup{\arabic*.},leftmargin=*]
\item If no component is large, introduce one variable for the sum of
 all component contributions and impose $e\ge(n-1)/14$.
\item If at least two components are large, select two distinct large
 components. Introduce a separate variable for each, with lower bound
 $1/2$, and a third variable for the sum over the remaining components.
 The exact ranks of the selected components need not be specified:
 their orthogonality and common lower bound suffice.
\item If exactly one component is large, of type $T$ and rank $r$, put
 $S=n-r$. Apply any local exclusion whose hypotheses hold. Otherwise,
 if $S\ge8$, use $e_T\ge1/2$ and
 $e_{\mathcal S}\ge(S-1)/14$.
\item If $0<S<8$, let $d(T,r)$ be the strongest applicable component
 bound in Section~\ref{c:large} that allows arbitrary other finite
 components. Estimates requiring a particular additional component
 are not used to define $d(T,r)$. If
 \begin{equation}\label{c:aggregate-test}
                  d(T,r)+(S-1)/14\ge1,
 \end{equation}
 use the class specified by these two bounds on the large component
 and the sum over the remaining components.
\item In every remaining case, list every spherical forest of rank
 $S$ with no large component, first without an exceptional configuration
 and then with each exceptional configuration that is allowed on that forest. Include all ordinary, large-component,
 joint and additional-component estimates that apply to this class,
 and the local exclusions of Section~\ref{c:local-exits} with their
 stated hypotheses. For $S=0$, the full vertex diagram consists of
 the one large component.
\end{enumerate}

The classes described by sums of components are included in the linear
program. Their displayed inequalities, each with multiplier one,
already imply a total contribution at least one; this includes
$S=8$ in step three. The test \eqref{c:aggregate-test} determines
when individual remaining components need to be listed. It changes
neither the hypotheses nor the validity of any local estimate.

\begin{proposition}[Exhaustiveness of the compact configurations]
\label{c:complete-cover}
Let a compact Coxeter polytope of dimension $18\le n\le29$ and a
farthest vertex be given. Its vertex configuration either satisfies
the hypotheses of one of the stated local exclusions or belongs to
one of the classes just constructed. The squared projection norms
of $x$ satisfy every inequality assigned to that class.
\end{proposition}
\begin{proof}
The spherical diagram has no large component, at least two large
components, or exactly one. The first two cases give steps one and
two. In the third case, its type and rank determine $S$.
Lemma~\ref{c:spectator} and the applicable component estimate justify
steps three and four when their numerical conditions hold.
Otherwise step five includes the remaining component multiset and
its exceptional configuration, if any. Components of equal type can
be relabelled, together with their facets, to place a marked component
first. This relabelling does not change any other parameter: each
ordinary dihedral component keeps its full integer range, and an
$O_7$ mark fixes the label seven only on its incident component.
Thus every inequality and local exclusion is applied with its required
component and facet hypotheses. Orthogonality gives
\eqref{c:budget} for the fixed vector $x$.
\end{proof}

For the explicit enumeration in step five, use the following ordered
list of component types through rank seven:
\[
\begin{array}{c|l}
\text{rank}&\text{types}\\\hline
1&A_1\\
2&A_2=I_2(3),\ I_2(m)\ (m\ge4)\\
3&A_3,B_3,H_3\\
4&A_4,B_4,D_4,F_4,H_4\\
5&A_5,B_5,D_5\\
6&A_6,B_6,D_6,E_6\\
7&B_7,D_7,E_7.
\end{array}
\]
Choose types in nondecreasing order until their ranks sum to $S$.
Then mark no exceptional configuration, one $A_1$ as $I$, two $A_1$
components as $J$, one $A_2$ as $N_7$, or one dihedral component as
$O_7$. Ordering and permutation of equal components avoid repetitions.
All dihedral parameters remain universally quantified. This gives the
$829$ classes used below. Their exhaustiveness follows from
Proposition~\ref{c:complete-cover}, independently of this count.

\subsection{Local exclusions}
\label{c:local-exits}

Table~\ref{c:local-exit-table} states the configuration classes
excluded by local theorems of Supplement~C. Their scaling coefficient
is zero in the normalized linear program; they are not assigned an
artificial lower bound on a component norm. Except for C1, each displayed forest is the full Coxeter diagram
at the chosen vertex, and the statements are used only for these full diagrams.
Every Coxeter label and parameter for separated facets ranges over
the full set in the cited local theorem.

\begin{longtable}{>{\raggedright\arraybackslash}p{.08\linewidth}
 >{\raggedright\arraybackslash}p{.43\linewidth}
 >{\raggedright\arraybackslash}p{.39\linewidth}}
\caption{Local exclusions in the compact case, proved in Supplement~C.}\label{c:local-exit-table}\\
\toprule
Code&Full vertex diagram or class&Local theorem\\\midrule
\endfirsthead
\toprule Code&Full vertex diagram or class&Local theorem\\\midrule\endhead
C1&Any $A_r$, $21\le r\le29$; $B_r$, $23\le r\le29$;
 or $D_r$, $24\le r\le29$, with arbitrary other finite components&
 Upper bound on the rank of a component\\
C2&$B_r+I$, $r\ge20$, or $D_r+I$, $r\ge21$&
 The classical-component theorem in the presence of $I$\\
C3&$A_{16}+N_7$ or $A_{16}+O_7$&
 The endpoint theorems for $N_7$ and $O_7$ in dimension $18$\\
C4&$A_{17}+A_1$ or $A_{17}+I$, with ordinary $A_1$&
 The theorem for a rank-one complement in dimension $18$\\
C5&$B_{16}+N_7/O_7$ or $D_{16}+N_7/O_7$&
 The rank-sixteen theorems for $N_7$ and $O_7$, including both signs at the $D$ fork\\
C6&$B_{17}+I$&The theorem for a short isolated opposite facet associated with the $A_1$ component\\
C7&$D_{17}+A_1$ or $D_{17}+I$, with ordinary $A_1$&
 The separate theorems for an ordinary $A_1$ and for $I$\\
C8&$A_{18}$&The endpoint exclusion for the indicated full diagram\\
C9&$B_{18}$&The terminal-facet theorem, including coincident and distinct opposite facets\\
C10&$D_{18}$&The endpoint and fork theorem, including intersecting, parallel and separated opposite facets\\
C11&$A_{17}+N_7/O_7$ or $A_{18}+I$&
 The separate local theorems in dimension nineteen\\
C12&$D_{17}+T_2$, $B_{17}+N_7/O_7$, or
 $B_{18}/D_{18}+A_1$, with ordinary $T_2,A_1$&
 The classical-component extension theorem in dimension nineteen\\
C13&$D_{17}+N_7/O_7$, $B_{18}+I$, or $D_{18}+I$&
 The theorems for $N_7$, $O_7$ and $I$ in dimension nineteen\\
C14&$A_{19}$, $B_{19}$ or $D_{19}$&
 The respective full-diagram theorems for types $A$, $B$ and $D$\\
C15&$A_{19}+I$&The theorem for $A_{19}$ with configuration $I$\\
C16&$B_{18}/D_{18}+N_7/O_7$ or $B_{19}/D_{19}+I$&
 The classical-component extension theorem in dimension twenty\\
C17&$A_{20}$, $B_{20}$ or $D_{20}$&
 The respective full-diagram theorems in dimension twenty\\
C18&$A_{19}+N_7/O_7$ or $A_{20}+I$&
 The $A$-component theorems for $N_7$, $O_7$ and $I$\\
C19&$D_{19}+N_7/O_7$ or $D_{20}+I$&
 The $D$-component extension theorem in dimension twenty-one\\
C20&$B_{21}$ or $D_{21}$&
 The respective full-diagram theorems in dimension twenty-one\\
C21&$B_{22}$, $D_{22}$ or $D_{23}$&
 The terminal-facet theorems for the indicated full diagrams\\\bottomrule
\end{longtable}

For example, C10 includes every terminal-facet and fork alternative
for the full diagram $D_{18}$. Auxiliary vertices used in a local
proof do not contribute additional terms to \eqref{c:budget}.

\subsection{The dual certificate}
\label{c:certificate}

The construction gives $829$ configuration classes, of which $151$
are ruled out by local exclusions. For each remaining class the dual
has one nonnegative multiplier for each applicable norm inequality.
The coefficient of every component contribution is at most one.
The rational solution satisfies
Theorem~\ref{thm:energy-certificate} simultaneously for all classes.
In particular, inequalities involving $I$ use the same variable $e_I$.

\begin{theorem}[The compact dimension bound]\label{c:recovery}
The local estimates of Sections~\ref{c:small}, \ref{c:large}
and~\ref{c:local-exits}, together with spherical classification,
make the normalized linear program infeasible in dimensions
$18,\ldots,29$. Together with Vinberg's upper bound $29$, they
imply $n\le17$ for compact hyperbolic Coxeter polytopes.
\end{theorem}
\begin{proof}
By Proposition~\ref{c:complete-cover}, a hypothetical polytope belongs
to one of the listed classes. A locally excluded class is impossible.
For any other class, Theorem~\ref{thm:energy-certificate} and the
rational dual solution imply that the total squared projection norm
is at least one, contrary to \eqref{c:budget}. The exact verification
checks all $3157$ dual variables and $2224$ inequalities; the records
are given in the reproduction appendix. Vinberg's theorem excludes
$n\ge30$.
\end{proof}

The configuration list is obtained from spherical classification and
the stated component estimates. It does not assume the earlier list
of eleven forests at rank eighteen.

\section{The compact arithmetic case}
\label{a:section}

For a compact arithmetic Coxeter polytope, the reflection group is
an arithmetic lattice. We apply the additional local results of
Supplement~A of~\cite{BoundsSupplement}, keeping their hypotheses
on the defining field and the normal lines of facets. The norm
inequality remains \eqref{c:budget}. The preceding compact classes
apply in dimensions $18,\ldots,29$; we add arithmetic classes at
ranks $16$ and $17$.

\subsection{Arithmetic hypotheses and local results}
\label{a:inputs}

In the quadratic-form model used in Supplement~A, $K$ is a totally
real number field, $q$ is a quadratic form over $K$, and $B$ is its
associated symmetric bilinear form. At one distinguished real
embedding $\sigma_0:K\hookrightarrow\R$, the form $q$ has Lorentzian
signature and defines the hyperbolic space containing the polytope.
At every other real embedding $\sigma$, the conjugate form
$q^\sigma$ is positive definite.

Choose vectors $\alpha_i$ on the facet-normal lines over this
same field; these are the root vectors used below. The associated
Cartan entries are
\begin{equation}\label{a:cartan}
 C_{ij}=\frac{2B(\alpha_i,\alpha_j)}{B(\alpha_i,\alpha_i)}.
\end{equation}
The local arithmetic theorems require the cyclic products
\begin{equation}\label{a:cycles}
             C_{i_1i_2}C_{i_2i_3}\cdots C_{i_\ell i_1}
\end{equation}
to belong to the ring of integers $\mathcal O_K$.
Replacing $\alpha_i$ by $t_i\alpha_i$, with $t_i\in K^\times$,
multiplies $C_{ij}$ by $t_j/t_i$; these factors cancel in each
cyclic product. The field reductions used below are those proved in
Supplement~A with these hypotheses.

At an embedding $\sigma\ne\sigma_0$, linearly independent conjugate
normal vectors have a positive-definite Gram matrix. The Gram matrix
of all facet normals may be redundant, and is then only positive
semidefinite. Normalizing each vector to unit length may introduce
square roots. The resulting field of real cosine entries therefore
does not replace the common field $K$ or the square classes of the
root norms. The norm class of $\alpha_i$ is the class of
$B(\alpha_i,\alpha_i)$ in $K^\times/(K^\times)^2$;
rescaling by $t_i\in K^\times$ multiplies this norm by $t_i^2$
and leaves its class unchanged.

Applications of Bugaenko's theorem~\cite{Bugaenko1984} use the
specified invariant lattice $L$ and the proved inclusion
$\Gamma(P)\subseteq\operatorname{Ref}(L)$. Here $\Gamma(P)$ is
the reflection group of the polytope, and $\operatorname{Ref}(L)$
is the full reflection group of that lattice. Equality of defining
fields alone does not establish the inclusion.

At ranks $16$ and $17$, Supplement~A gives necessary lists of
configurations after its local arithmetic reductions. We take those
lists as proved local results. The linear program combines their
norm inequalities and local exclusions; it does not derive the
preceding reductions.

\begin{assumption}[Local arithmetic results]
\label{a:input-package}
The following results of Supplement~A are used with all of their
stated hypotheses.
\begin{enumerate}[label=\textup{(\roman*)},leftmargin=*]
\item The common totally real field, facet-normal lines, integral
 Cartan cycles, root-norm square classes and positivity at every
 real embedding other than the distinguished one are as required by each local theorem.
\item The fourteen rank-seventeen families and nineteen rank-sixteen
 families in Tables~\ref{a:rank17} and~\ref{a:rank16} are exhaustive
 necessary lists, including all exceptional configurations, other
 components and field alternatives occurring in their proofs.
\item The estimates in Section~\ref{a:costs} and
 Tables~\ref{a:rank17}--\ref{a:rank16} hold for the indicated
 components of the full vertex diagrams. Each component bound includes
 every contribution used in its local proof.
\item The local exclusions of Table~\ref{a:local-exits} hold for the
 stated facet configurations, including boundary values and the
 identification of each opposite facet associated with an $A_1$ component.
\item The four finite classifications in
 Table~\ref{a:computed-inputs} are exhaustive on their stated
 parameter sets. The geometric reductions prove that the relevant
 facet configurations belong to those sets.
\item Each application of the golden-lattice theorem uses the lattice
 and full reflection-group inclusion specified above.
\end{enumerate}
\end{assumption}

These statements contain arithmetic information in addition to
positive semidefiniteness of real Gram matrices.

\begin{table}[htbp]
\centering
\caption{Four finite classifications used in Supplement~A.
An algebraic solution need not correspond to a polytope.}
\label{a:computed-inputs}
\begin{tabular}{>{\raggedright\arraybackslash}p{.08\linewidth}
 >{\raggedright\arraybackslash}p{.43\linewidth}
 >{\raggedright\arraybackslash}p{.39\linewidth}}
\toprule
No.&Classified parameter set&Conclusion used in the local proof\\\midrule
(i)&The quadratic-field parameter set for full diagram $A_{17}$&
 No solutions\\
(ii)&The two cubic-field parameter sets for full diagram $A_{17}$,
 each with the required endpoint contact at a finite Coxeter angle&
 No solutions in either set\\
(iii)&The common quadratic parameter set obtained by separate
 reductions for full diagrams $B_{17}$ and $D_{17}$&Two possible
 pairs of parameters; each is then subject to the appropriate
 $B$- or $D$-component norm estimate\\
(iv)&The quadratic parameter set obtained from full diagram $B_{16}$&
 One possible pair of parameters, to which the estimate along the
 specified edge is applied\\\bottomrule
\end{tabular}
\end{table}

In case (ii), the fields are $\Q(2\cos(2\pi/7))$ and
$\Q(2\cos(2\pi/9))$. The required finite contact has Coxeter
label seven or eighteen, respectively. Thus an endpoint opposite
wall meets the corresponding removed incident wall at angle
$\pi/7$ or $\pi/18$, respectively. These
are the two cases of the $A_{17}$ endpoint reduction in Supplement~A.

The exclusion of the full diagram $D_{16}$ is proved separately by
a direct argument over dyadic completions. It does not use one of
the four finite classifications in Table~\ref{a:computed-inputs}.
It is distinct from the exclusion of $D_{16}+I$ at rank seventeen.

\subsection{Notation and local norm inequalities}
\label{a:costs}

The quantities $e_C$ continue to denote squared norms of orthogonal
projections of the fixed vector $x$. Write $T_2$ for an ordinary
irreducible $I_2(m)$, $m\ge3$, $T_3$ for $A_3,B_3$ or $H_3$,
and $T_4$ for one of $A_4,B_4,D_4,F_4,H_4$.
Let $\mathcal K$ denote either exceptional configuration $N_7$ or
$O_7$; they are kept as separate alternatives. An unqualified $A_1$
is ordinary and nonshort. The configurations $I,J,N_7,O_7$ are those
of Section~\ref{c:small}.

Set
\begin{align*}
\mathcal R_5={}&\{A_5,D_5,A_4+A_1,H_4+A_1,
                   T_3+T_2,2T_2+A_1\},\\
\mathcal R'_5={}&\{A_5,D_5,A_4+A_1,H_4+A_1,
                   A_3+T_2,2T_2+A_1\}.
\end{align*}
Every member contributes at least $5/14$ by
Table~\ref{c:ordinary-table}. The restriction to $A_3$ in
$\mathcal R'_5$ is part of the rank-sixteen classification.

In addition to the ordinary estimates of
Table~\ref{c:ordinary-table}, Supplement~A gives
\begin{align}
 e_{I+T_2}&>21/100,& e_{I+T_3}&>2/7,&
 e_{I+T_4}&>23/70,\label{a:I-joints}\\
 e_{I+A_1}&>5051/25000=.20204.\label{a:I-A1}
\end{align}
It also gives $e_{T_2}>1/5$ for an ordinary component outside
$N_7$ or $O_7$, and $e_{A_1}>2/15$ for an ordinary component
outside $O_7$.

Tables~\ref{a:rank17} and~\ref{a:rank16} state the remaining local
bounds on the full vertex diagrams in their middle columns. Every
listed dihedral parameter and field alternative ranges over the set
specified in the local theorem. Separate summands refer to estimates
on disjoint component sets, unless a joint estimate is explicitly
named. An estimate assigned to the large component bounds that
component alone, even when its proof uses other facets. For example,
$16/19$ for $A_{14}$ with an ordinary $T_2$ and
$.91286445198$ with $N_7$ or $O_7$ are different statements.
The ordinary $T_2$ estimate $3/19$ follows from the arithmetic
dihedral theorem when $m\ge4$, and from $319/2000>3/19$ when
$m=3$. For the indicated $B/D$ components with $N_7$ or $O_7$,
the bounds $.91859$ and $.93859$ already include the correction
$.0487$ from the opposite facets. That correction is part of the
bound on the same component, rather than a separate summand.
We display the inequalities used by the certificate below. Additional
valid inequalities in the electronic model may be assigned multiplier
zero.

\subsection{The rank-seventeen configurations}
\label{a:rank17-section}

Table~\ref{a:rank17} gives the fourteen families from Supplement~A.
Braces and semicolons separate alternatives, all of which are included.
The local exclusions X1--X13 are stated in Table~\ref{a:local-exits}.

\begin{longtable}{>{\raggedright\arraybackslash}p{.06\linewidth}
 >{\raggedright\arraybackslash}p{.35\linewidth}
 >{\raggedright\arraybackslash}p{.49\linewidth}}
\caption{The rank-seventeen families and the inequalities used in the dual solution.}
\label{a:rank17}\\
\toprule
No.&Full vertex diagram&Norm bounds with sum at least one, or local exclusion\\\midrule
\endfirsthead
\toprule No.&Full vertex diagram&Norm inequalities or local exclusion\\\midrule\endhead
1&$A_{11}+I+(A_5\text{ or }D_5)$&
 $13/20+5/14=141/140$; the nonnegative contribution from $I$ is omitted\\
2&$A_{11}+\mathcal K+U_4$, $U_4\in\{A_4,B_4,D_4,F_4\}$&
 $13/20+9/100+2/7$\\
3&$A_{12}+R$, $R\in\mathcal R_5$&
 $5/7+5/14$\\
4&$A_{12}$ plus $I+T_4$; $I+T_3+A_1$; or $I+2T_2$&
 $5/7$ plus, respectively, the joint $23/70$;
 the joint $2/7$ and $1/10$; or the joint $21/100$ and $1/7$\\
5&$A_{14}$ plus $T_3$; $T_2+A_1$; or $I+T_2$&
 $4/5$ plus, respectively, $1/4$; $1/7+1/10$;
 or the joint $21/100$\\
6&$A_{15}$ plus $T_2$; $2A_1$; $I+A_1$; $N_7$; or $O_7$&
 $T_2$: X1. The next two give $.8348+.2$ and
 $.8348+.20204$. Either exceptional configuration gives $.9191408+.09$\\
7&$A_{16}$ plus ordinary $A_1$ or configuration $I$&X2 for both rank-one alternatives\\
8&$A_{17}$&X3, using classifications (i), (ii) and the reductions to their parameter sets\\
9&$B_{15}+N_7$ or $B_{15}+O_7$&
 $.93859+.09=1.02859$\\
10&$B_{16}+I$&X4\\
11&$B_{17}$&$e_{B_{17}}\ge1.00184$, after the
 separate reduction and norm estimates for both parameter pairs\\
12&$D_{15}+N_7$ or $D_{15}+O_7$&
 $.91859+.09=1.00859$\\
13&$D_{16}$ plus ordinary $A_1$ or configuration $I$&
 X5 for ordinary $A_1$; X6 for $I$\\
14&$D_{17}$&$e_{D_{17}}\ge1.00184$, after its
 reduction for the $D$ fork and its norm estimates for both parameter pairs\\\bottomrule
\end{longtable}

The fourteen families expand to $39$ parametrized configuration classes.
Some may not be geometrically realizable. Repeated symbols such as
$2T_2$ denote components with independent parameters. For $I+2T_2$,
apply the joint estimate to $I$ and one dihedral component, and the
ordinary estimate to the other dihedral component. Thus $e_I$ occurs
with coefficient one.

\subsection{The rank-sixteen configurations}
\label{a:rank16-section}

Table~\ref{a:rank16} gives the nineteen rank-sixteen families,
including the restricted family $\mathcal R'_5$.

\begin{longtable}{>{\raggedright\arraybackslash}p{.06\linewidth}
 >{\raggedright\arraybackslash}p{.35\linewidth}
 >{\raggedright\arraybackslash}p{.49\linewidth}}
\caption{The rank-sixteen families and the inequalities used in the dual solution.}
\label{a:rank16}\\
\toprule
No.&Full vertex diagram&Norm bounds with sum at least one, or local exclusion\\\midrule
\endfirsthead
\toprule No.&Full vertex diagram&Norm inequalities or local exclusion\\\midrule\endhead
1&$A_{10}+I+(A_5\text{ or }D_5)$&$404/625+5/14$\\
2&$A_{10}+\mathcal K+U_4$, $U_4\in\{A_4,B_4,D_4,F_4\}$&
 $404/625+9/100+2/7$\\
3&$A_{11}+R$, $R\in\mathcal R'_5$&$13/20+5/14=141/140$\\
4&$A_{11}+I+T_4$&The joint $A_{11}+I$ bound $18/25$,
 plus the disjoint $T_4$ bound $2/7$, giving $176/175$\\
5&$A_{11}+O_7+A_3$&$13/20+9/100+13/50=1$\\
6&$A_{12}$ plus $A_4$; $H_4$; or $2T_2$&$5/7+2/7=1$\\
7&$A_{12}+I+T_3$&$5/7$ plus the joint $I+T_3$ bound $2/7$\\
8&$A_{12}+\mathcal K+T_2$&$5/7+9/100+1/5$, using the stronger ordinary-component bound in the presence of $N_7$ or $O_7$\\
9&$A_{13}$ plus $A_3$; $T_2+A_1$; $I+T_2$; or $O_7+A_1$&
 Respectively: $53/70+1/4$; $53/70+1/7+1/10=1$;
 X7; or $233/300+9/100+2/15=1$\\
10&$A_{14}$ plus $T_2$; $2A_1$; $I+A_1$; $N_7$; or $O_7$&
 Respectively: $16/19+3/19=1$; $4/5+.2=1$;
 $4/5+.20204$; or $.91286445198+.09$ with either $N_7$ or $O_7$\\
11&$A_{15}$ plus ordinary $A_1$ or configuration $I$&X8 for every isolated and nonisolated rank-one alternative\\
12&$A_{16}$&X9\\
13&$B_{14}+N_7$ or $B_{14}+O_7$&
 $.91859+.09=1.00859$\\
14&$B_{15}$ plus ordinary $A_1$ or configuration $I$&
 X10 for ordinary $A_1$; X11 for $I$\\
15&$B_{16}$&$e_{B_{16}}\ge1.00640742675$, after
 its reduction to the parameter set in (iv) and the edge estimate for the single solution\\
16&$D_{13}$ plus $I+T_2$ or $O_7+A_1$&
 $.82436$ plus the bound $.21$ on the remaining components; joint on $I+T_2$,
 and coarsened from $.09+2/15$ on $O_7+A_1$\\
17&$D_{14}$ plus ordinary $I_2(m)$, $m\ge4$; $N_7$; or $O_7$&
 Ordinary: $.84816+3/19$. Either exceptional configuration:
 $.91859+.09$\\
18&$D_{15}$ plus ordinary $A_1$ or configuration $I$&
 X10 for ordinary $A_1$; X12 for $I$\\
19&$D_{16}$&X13, the quadratic-field reduction and dyadic theorem\\\bottomrule
\end{longtable}

Every displayed sum is at least one and therefore contradicts
\eqref{c:budget}. No uniform positive difference from one is needed.
For example, family 5 uses three estimates on disjoint component sets
whose bounds sum to one. Family 4 uses the joint estimate on
$A_{11}+I$:
\[
 e_{A_{11}}+e_I\ge18/25,\qquad e_{T_4}\ge2/7.
\]
The component sets in these two inequalities are disjoint.
The bound $18/25$ applies to the sum on $A_{11}+I$; it cannot be
assigned to $A_{11}$ alone or added to another estimate with a second
full coefficient on $e_I$.

The nineteen families give $49$ parametrized configuration classes. Together
the two lists contain $88$ alternatives: $72$ are treated by norm
inequalities, and $16$ by local exclusions. The dihedral parameters
retain their full, possibly unbounded ranges.

\subsection{Local arithmetic exclusions}
\label{a:local-exit-section}

Table~\ref{a:local-exits} states every local exclusion used in
Tables~\ref{a:rank17} and~\ref{a:rank16}. Each is a local theorem
of Supplement~A with the stated arithmetic and geometric hypotheses.
It sets the scaling coefficient of that configuration class to zero in the normalized
linear program, rather than an additional norm estimate.

\begin{longtable}{>{\raggedright\arraybackslash}p{.07\linewidth}
 >{\raggedright\arraybackslash}p{.31\linewidth}
 >{\raggedright\arraybackslash}p{.52\linewidth}}
\caption{Local arithmetic exclusions. The indicated diagram is the full vertex diagram.}\label{a:local-exits}\\
\toprule
Code&Full vertex diagram&Local theorem\\\midrule
\endfirsthead
\toprule Code&Full vertex diagram&Local theorem\\\midrule\endhead
X1&$A_{15}+T_2$, rank $17$&The ordinary-dihedral theorem, including all defining-field and small-norm alternatives.\\
X2&$A_{16}+A_1$ or $A_{16}+I$, rank $17$&The theorem for opposite facets associated with the $A_1$ component, including isolated and nonisolated cases.\\
X3&$A_{17}$, rank $17$&The endpoint and field reductions to the finite parameter sets in (i) and (ii).\\
X4&$B_{16}+I$, rank $17$&The theorem for a short isolated opposite facet, for both quadratic and cubic fields.\\
X5&$D_{16}+A_1$, rank $17$, with ordinary $A_1$&The ordinary rank-one theorem, including both alternatives at the fork.\\
X6&$D_{16}+I$, rank $17$&The theorem for this full diagram over fields of degrees two, three and four, including both fork signs and every dyadic case.\\
X7&$A_{13}+I+T_2$, rank $16$&The endpoint theorem involving the remaining components, with every ordinary dihedral label.\\
X8&$A_{15}+A_1$ or $A_{15}+I$, rank $16$&The rank-one opposite-facet theorem, including every height boundary.\\
X9&$A_{16}$, rank $16$&The endpoint and field theorem for the full diagram, including its finite-contact cases and ratios of Lorentz time coordinates.\\
X10&$B_{15}+A_1$ or $D_{15}+A_1$, rank $16$, with ordinary $A_1$&
 The respective nonisolated reductions and ordinary isolated theorem, including $h=1/2$.\\
X11&$B_{15}+I$, rank $16$&The isolated-short theorem
 for $0<h<1/2$.\\
X12&$D_{15}+I$, rank $16$&The isolated-short theorem
 for $0<h<1/2$, including both signed forks.\\
X13&$D_{16}$, rank $16$&The quadratic-field reduction and the dyadic theorem, including the odd, inert, ramified and split cases.\\\bottomrule
\end{longtable}

For each rank-one theorem, fix the opposite facet associated with
the indicated $A_1$ component. If its spatial normal has a nonzero
projection to the large component, the nonisolated case applies.
Otherwise the facet is isolated, and the parameter given by the square
of the coefficient of $p$ in its unit normal satisfies $0<h<1/2$ or
$h\ge1/2$. The equality case is included in the ordinary theorem.
Thus these alternatives concern facets of the polytope, not arbitrary
formal extensions of its diagram.

For the full diagrams $B_{17}$ and $D_{17}$, the local reductions
separately place the parameters in class (iii) of
Table~\ref{a:computed-inputs}. The component estimates are then
applied to both possible parameter pairs. For $B_{16}$, the reduction
to class (iv) and the estimate along its specified edge treat the
single possible pair. These finite classifications do not give empty
parameter sets. The full diagram $D_{16}$ is instead excluded by the
dyadic theorem X13.

\subsection{The dual certificate for the arithmetic case}
\label{a:certificate}

The arithmetic tables give an explicit dual solution. For every
alternative treated by norm inequalities, assign multiplier one to
the displayed inequalities on disjoint component sets and zero to
any additional inequalities. Each component variable then has total
coefficient at most one, while the sum of the lower bounds is at
least one. Alternatives X1--X13 have zero scaling coefficient and
therefore zero normalized deficit, by their local exclusion theorems.

The computation uses the same linear-program construction for all
$88$ arithmetic alternatives and the $829$ compact classes in
dimensions $18$--$29$. There are $917$ classes in total, of which
$167$ are locally excluded. The checked dual has $3324$ variables
and $2488$ inequalities. The reproduction appendix records these
linear verifications separately from the local geometric proofs.

\begin{theorem}[The compact arithmetic dimension bound]\label{a:recovery}
Assume the local results of Input~\ref{a:input-package}, the compact
local results of Section~\ref{c:section}, and Vinberg's compact upper
bound $29$. The combined normalized linear program for dimensions
$16,\ldots,29$ is infeasible. Consequently a compact arithmetic
hyperbolic Coxeter polytope satisfies $n\le15$.
\end{theorem}
\begin{proof}
At ranks $16$ and $17$, every arithmetic facet configuration belongs
to one of the respective necessary lists. A locally excluded class
is impossible. On each remaining class, the dual inequalities force
the sum of squared projection norms to be at least one, contrary to
\eqref{c:budget}. In dimensions $18$--$29$, the compact classes
and Theorem~\ref{c:recovery} apply without change. Vinberg's theorem
excludes dimensions $30$ and above.
\end{proof}

\section{Exclusion in all higher dimensions}
\label{sec:thresholds}

Proposition~\ref{prop:threshold} applies to each of the three fixed
families of configurations. The following table summarizes the
additional estimates needed to cover the full dimension range.

\begin{center}\small
\begin{tabular}{>{\raggedright\arraybackslash}p{.17\linewidth}
 >{\raggedright\arraybackslash}p{.35\linewidth}
 >{\raggedright\arraybackslash}p{.38\linewidth}}
\toprule
Case&Finite calculation&Higher dimensions\\\midrule
Finite volume&Finite and ideal vertices at $n=58$&The analytic face
 and vertex estimates, together with the stated reduction of ideal
 links, apply at every $n\ge59$\\
Compact&Exhaustive classification at $18\le n\le29$&Vinberg's
 theorem gives $n\le29$\\
Compact arithmetic&Arithmetic classifications at $n=16,17$, and
 the compact classes at $18\le n\le29$&The same compact theorem
 gives $n\le29$\\\bottomrule
\end{tabular}
\end{center}

Vinberg's bound excludes compact polytopes in dimensions
$n\ge30$~\cite{Vinberg1984,Vinberg1985}. The arithmetic program
uses the compact classes and inequalities at dimensions $18$--$29$
unchanged. The general finite-volume case instead uses the analytic
estimates proved for every higher dimension. Thus the certificates
exclude all dimensions at least $58$, $18$ and $16$, respectively.
This conclusion requires the stated higher-dimensional estimates;
infeasibility of an isolated exact-dimension problem would not by
itself imply it.

A lower threshold would require an exhaustive family of classes
and valid local inequalities in the newly included dimensions.
The monotonicity proposition does not assert that an arbitrary
face of a Coxeter polytope is again a Coxeter polytope.

\section{Exact verification of the linear certificates}
\label{sec:certificates}

For the compact cases, the matrix records the support, rational lower
bound and applicability class of each local inequality. The dual
programs for all classes are solved together. A numerical solver
provides multipliers, which are converted to rational numbers and
checked against every matrix row and variable bound by exact
arithmetic. A separate run regenerates the configuration lists and
both the primal and dual matrices, compares the resulting files
byte for byte, and checks the rational solution again.

\begin{center}\small
\begin{tabular}{lrrrr}
\toprule
Problem&Classes&\shortstack{Local\\exclusions}&
\shortstack{Primal\\variables/rows}&\shortstack{Dual\\variables/rows}\\
\midrule
Compact&829&151&$3471/4876$&$3157/2224$\\
\shortstack[l]{Arithmetic and\\compact classes}&917&167&$3869/5147$&$3324/2488$\\
\bottomrule
\end{tabular}
\end{center}

In finite volume the local matrix has $5\,810$ variables and
$103\,606$ rows. The rank-$58$ forest inequalities are represented
by the complete collection of finite networks. The local matrix
and every network use the same weight coefficients. Their exact
feasible assignment gives the upper bound $377$, and the eight-row
dual combination gives the matching lower bound. These checks
used no numerical optimizer and kept the network inequalities in
separate matrices. The subsequent verification repeated the local
inequalities, the completeness checks joining the networks, and
the lower-bound certificate, using the previously checked networks.

An inequality may be redundant on a class. For example, adding a
positive correction to the initial-segment norm estimate gives a
stronger inequality on the same support. Either may remain in the
matrix. Conversely, a zero multiplier in one dual solution proves
only that this solution does not use the corresponding inequality;
it does not prove redundancy of that inequality.

\section{Further constraints}
\label{sec:strengthening}

Additional necessary conditions may restrict a configuration class
or add inequalities to its cone. Additional functions of a corner
may enlarge the class of weight functions available in finite
volume. These modifications preserve the argument when every
polytope still determines a feasible point and the configuration
data remain fixed as the dimension threshold varies.

Some local results could instead be expressed by semidefinite
constraints on Gram matrices. Such a formulation would have to
preserve the geometric applicability of every constraint. In the
arithmetic case it would also have to account for integral Cartan
cycles, conjugate forms at the other real embeddings, and the
specified reflection-group inclusion in a lattice. Real Gram
positivity alone does not imply those arithmetic conditions.

\appendix
\section{Sources and reproducibility}
\label{app:reproduction}

The geometric results are given in the dimension-bound article and
its supplement~\cite{BoundsArticle,BoundsSupplement}. The
finite-volume weight function also uses the comparisons between
its two dimension ranges and the analytic estimates described in
the research notes. The repository~\cite{LPRepository} provides
the reading guide \path{notes/lp/UNIFORM_LP_PARADIGM.md} and the
result summary
\path{notes/lp/existing_proof_lp_recovery_results_20261118.md}.
The principal source files and verification records are listed
below. The repository's reproduction guide distinguishes the
publicly available files from the larger research data not yet
distributed with the source checkout.

\begingroup\small\raggedright
\begin{longtable}{>{\raggedright\arraybackslash}p{.23\linewidth}
 >{\raggedright\arraybackslash}p{.70\linewidth}}
\toprule
Record&Path\\\midrule\endfirsthead
\toprule Record&Path\\\midrule\endhead
Finite-volume verification&\path{audit/joint58_complete_lp_execution_20261110.json}\\
Subsequent verification&\path{audit/lp_rerun_report_execution_20261115.json}\\
Compact configuration generator&\path{scripts/compact17_proof_recovery_energy_catalogue_20261118.py}\\
Arithmetic configuration generator&\path{scripts/existing_proof_arithmetic_energy_domains_20261118.py}\\
Common LP construction&\path{scripts/run_existing_proof_energy_lp_20261118.py}\\
Exact matrix checker&\path{scripts/finite_volume_lp_matrix_20261106.py}\\
Compact matrices and certificates&\path{results/existing_proof_compact17_energy_master01_20261118/}\\
Arithmetic matrices and certificates&\path{results/existing_proof_arithmetic15_energy_master01_20261118/}\\
Independent verification review&\path{audit/existing_proof_energy_lp_execution_independent_review_20261118.json}\\
Verification summary&\path{audit/existing_proof_lp_recovery_completion_root_acceptance_20261118.json}\\
\bottomrule
\end{longtable}
\endgroup

Each compact output directory contains the configuration list,
normalized primal problem, dual problem, rational solution, and
construction and verification reports. The execution records
include the commands, source hashes, resource limits and measurements.
The mathematical applicability of each local inequality is checked
in the source review; the rational checker verifies its linear
combination. File hashes identify the exact versions used.

Reproduction uses the repository's serial monitored runner. The
\texttt{export-solve} stage constructs the problem and finds a dual
solution; a separate \texttt{verify} stage regenerates both matrices
and checks the rational certificate. For arithmetic, both the
arithmetic and compact configuration generators are supplied.
Each run requires fresh output and monitor paths. The execution
reports give the exact commands and resource limits.

The resulting certificates prove the linear consequences of the
stated local theorems and classifications. The geometric proofs,
field reductions, four finite arithmetic classifications and lattice
theorems retain the hypotheses stated in the text. Their earlier
verification records are cited separately from the new global LP
checks.

\paragraph{Acknowledgments.}
The author used OpenAI's Codex as a research and coding assistant to
investigate proof ideas, develop the linear programming formulations,
improve the exposition, and generate and audit the accompanying Python
code. The author reviewed the resulting mathematics and code and takes
responsibility for the final content.

\bibliographystyle{unsrt}
\bibliography{references}
\end{document}